\documentclass[preprint,11pt]{elsarticle}

\usepackage{graphicx}
\usepackage[utf8]{inputenc}
\usepackage{multirow}
\usepackage{comment}
\usepackage{mdwlist}
\usepackage{amsfonts}
\usepackage{ulem}
\usepackage{amsmath} 
\usepackage{amsthm}
\usepackage{amssymb}
\usepackage{fancyhdr}
\newtheorem{lemma}{Lemma}
\newtheorem{theorem}{Theorem}
\usepackage{indentfirst}
\usepackage{booktabs}
\usepackage{verbatim}
\usepackage{xcolor}
\usepackage{latexsym}
\usepackage{comment}
\usepackage{amscd}
\usepackage{esvect}
\usepackage{graphicx}
\usepackage{upgreek}
\usepackage{caption}   
\usepackage[a4paper]{geometry}
\usepackage{mathrsfs}
\usepackage[numbers]{natbib}
\usepackage{float}
\usepackage{lipsum}
\definecolor{RoyalBlue}{RGB}{65, 105, 225}
\usepackage[export]{adjustbox} 
\usepackage{subfiles}
\usepackage{subfigure}
\usepackage{bm}
\usepackage{epstopdf}
\newcommand{\e}{\varepsilon}
\newcommand{\bx}{\boldsymbol{x}}

\newtheorem*{remark}{Remark}
\begin{document}
\begin{frontmatter}
\author{Jin Woo Jang$^{a}$}
\ead{jangjw@postech.ac.kr}
\author{Jae Yong Lee$^{b}$}
\ead{jaeyong@cau.ac.kr}
\author{Liu Liu$^{c}$}
\ead{lliu@math.cuhk.edu.hk}
\author{Zhenyi Zhu$^{c}$}
\ead{zyzhu@math.cuhk.edu.hk}
\address{$^{a}$Department of Mathematics, POSTECH (Pohang University of Science and Technology)}
\address{$^{b}$Department of AI, Chung-Ang University}
\address{$^{c}$Department of Mathematics, The Chinese University of Hong Kong}
\title{Deep learning-based moment closure for multi-phase computation of semiclassical limit of the Schr\"odinger equation}

\begin{abstract}
We present a deep learning approach for computing multi-phase solutions to the semiclassical limit of the Schrödinger equation. Traditional methods require deriving a multi-phase ansatz to close the moment system of the Liouville equation, a process that is often computationally intensive and impractical. Our method offers an efficient alternative by introducing a novel two-stage neural network framework to close the $2N\times 2N$ moment system, where $N$ represents the number of phases in the solution ansatz. In the first stage, we train neural networks to learn the mapping between higher-order moments and lower-order moments (along with their derivatives). The second stage incorporates physics-informed neural networks (PINNs), where we substitute the learned higher-order moments to systematically close the system. We provide theoretical guarantees for the convergence of both the loss functions and the neural network approximations. Numerical experiments demonstrate the effectiveness of our method for one- and two-dimensional problems with various phase numbers $N$ in the multi-phase solutions. The results confirm the accuracy and computational efficiency of the proposed approach compared to conventional techniques.
\end{abstract}
\end{frontmatter}

\allowdisplaybreaks

\section{Introduction}\label{sec.1}
We study the semiclassical limit of the Schrödinger equation with high-frequency initial data given by the Wentzel–Kramers–Brillouin (WKB) method
\begin{equation}
\label{Schro-eqn}
\begin{cases}
i \e \partial_t \psi^{\e}=-\frac{\e^2}{2} \Delta_x \psi^{\e} + \Phi(x) \psi^{\e},  \\[6pt]
\psi(0,x)=A_0(x) e^{ i S_0(x)/\e}, 
\end{cases}
\end{equation}
where $\psi^{\e}(t,x)$ is the wave function, $\e$ is the scaled Planck constant, and $\Phi(x)$ is the smooth potential function. From the wave function $\psi^{\e}$, one can derive fundamental physical observables, including the position density
\begin{equation}\label{density1} 
\rho(t,x) = | \psi^{\e}(t, x) |^2, 
\end{equation}
and the current density
\begin{equation}\label{density2}
J(t,x) = \e \operatorname{Im}\left(\overline{\psi^{\e}(t,x)}\nabla\psi^{\e}(t,x)\right) 
= \frac{1}{2i} \left( \overline{\psi^{\e}}\nabla \psi^{\e} - \psi^{\e}\nabla\overline{\psi^{\e}}\right). 
\end{equation}

Numerically solving the Schr\"odinger equation in the semiclassical regime presents several approaches. The Whitham averaging method \cite{Whitham1} derives modulation equations for multi-phase solutions, while the classical WKB approximation yields (to leading order) an eikonal equation for the phase and a transport equation for the amplitude \cite{Whitham2}. However, this Hamilton-Jacobi-type eikonal equation can develop singularities from smooth initial data, and viscosity-based numerical solvers may fail to capture the correct semiclassical limit beyond these singular points \cite{Engquist,Lions2}. Alternative methods like ray tracing have also been explored (see \cite{Ben1,Ben2}).

The Wigner transform \cite{GMMP,Lions,PR} offers a more robust framework, converting the problem into a linear Vlasov/Liouville equation in phase space that naturally handles caustics and provides globally valid multi-phase solutions. Kinetic schemes for this formulation have successfully validated such solutions numerically \cite{Jin-Xiantao}. Recent advances leverage machine learning to address moment closure problems, where neural networks effectively learn the implicit relationship between higher- and lower-order moments, a task often analytically intractable. This data-driven approach circumvents the need for explicit closure assumptions.

Regarding moment closure strategies, numerous models have been developed, including the classical $P_N$ model \cite{chandrasekhar1946radiative} and entropy-based $M_N$ model \cite{hauck2011high}. Recent advances in machine learning have introduced innovative approaches to this long-standing challenge:

\begin{itemize}
    \item Han et al. \cite{han2019uniformly} developed an autoencoder framework to extract optimal generalized moments for kinetic problems.
    \item Bois et al. \cite{MR4389617} proposed a CNN-based nonlocal closure for the Vlasov-Poisson system.
    \item Huang et al. \cite{huang2022machine,huang2023machine1,HLQW} established ML-based closures for radiative transfer and Boltzmann equations, preserving hyperbolicity through gradient learning.
\end{itemize}

Despite these advances, the application of machine learning to close the moment system arising from the Liouville equation (the semiclassical limit of the Schr\"odinger equation) remains unexplored. Our work presents the first data-driven closure framework specifically designed for this fundamental quantum mechanical system, offering new capabilities for semiclassical simulations.  

Our main contribution is a novel \textit{two-stage} deep learning framework for computing multi-phase solutions to the semiclassical limit of the Schrödinger equation. The key contributions of this work are as follows:
\begin{enumerate}
    \item We present an original two-stage physics-informed neural network (PINN) approach \cite{PINN1} for the multi-phase computations of the semiclassical limit of the  Schrödinger equation, which effectively solves the moment closure and  learns the intricate relationship between lower-order and highest-order moments;
    \item We provide rigorous theoretical guarantees for the convergence of both the loss function and the neural network approximations;
    \item We conduct comprehensive numerical validation that demonstrates the accuracy and efficiency of our method. 
\end{enumerate} This work establishes the first machine learning-based closure for the moment system of the Liouville equation in this context.

The rest of the paper is organized as follows. The rest of Section \ref{sec.1} derives the semiclassical limit through the Wigner transform and presents the motivation for our work. Section \ref{sec.2} develops our \textit{two-stage moment closure model}, detailing the neural network architecture, loss functions, and convergence analysis. Section \ref{sec.3} provides numerical experiments validating our method's ability to accurately close the moment system and compute multi-phase solutions. We conclude with discussion and future research directions in Section \ref{sec.4}. 

\subsection{The Wigner Transform and Semiclassical Limit}

For $f, g \in L^2(\mathbb{R}^d)$, we define the Wigner transform $W^\varepsilon(f,g)$ as the phase-space function:
\begin{equation}
W^\varepsilon(f,g)(t,x,v) = \frac{1}{(2\pi)^d} \int_{\mathbb{R}^d} \overline{f}\left(x+\frac{\varepsilon}{2}\sigma\right) g\left(x-\frac{\varepsilon}{2}\sigma\right) e^{i\sigma\cdot v} d\sigma,
\end{equation}
where $(t,x,v) \in [0,T] \times \Omega \times \mathbb{R}^d$ with $\Omega \subset \mathbb{R}^d$, and $\overline{f}$ denotes the complex conjugate. When applied to the solution $\psi^\varepsilon$ of the Schrödinger equation \eqref{Schro-eqn}, the Wigner function $W^\varepsilon(t,x,v) := W^\varepsilon(\psi^\varepsilon,\psi^\varepsilon)$ satisfies the Wigner equation \cite{BaoJin}:
\begin{equation}
\partial_t W^\varepsilon + v \cdot \nabla_x W^\varepsilon + \Theta[\Phi]W^\varepsilon = 0,
\end{equation}
where $\Theta[\Phi]$ is the pseudo-differential operator \cite{Lions}:
\begin{equation}
\Theta[\Phi]W^\varepsilon(t,x,v) = \frac{i}{(2\pi)^d} \int_{\mathbb{R}^{2d}} \delta\Phi^\varepsilon(x,\alpha) W^\varepsilon(t,x,\zeta) e^{-i\alpha\cdot(v-\zeta)} d\alpha d\zeta,
\end{equation}
with 
\begin{equation}
\delta\Phi^\varepsilon(x,\alpha) = \frac{1}{\varepsilon}\left[\Phi\left(x+\frac{\varepsilon}{2}\alpha\right) - \Phi\left(x-\frac{\varepsilon}{2}\alpha\right)\right].
\end{equation}
Under appropriate regularity conditions on $\Phi$ \cite{Sparber}, we have the semiclassical limit:
\begin{equation}
\delta\Phi^\varepsilon \xrightarrow{\varepsilon\to 0} \alpha \cdot \nabla_x \Phi(x).
\end{equation}
By Weyl's calculus \cite{Wigner1}, the limiting Wigner measure $w(t,x,v) := \lim_{\varepsilon\to 0} W^\varepsilon(\psi^\varepsilon,\psi^\varepsilon)$ satisfies the Vlasov equation:
\begin{equation}\label{Vlasov-eqn}
\partial_t w + v \cdot \nabla_x w - \nabla_x \Phi \cdot \nabla_v w = 0.
\end{equation}
For the WKB initial data $\psi(0,x) = A_0(x)e^{iS_0(x)/\varepsilon}$, the initial Wigner transform converges weakly to:
\begin{equation}
w(0,x,v) = \rho_0(x)\delta(v - u_0(x)),
\end{equation}
where $\rho_0(x) = |A_0(x)|^2$ and $u_0(x) = \nabla_x S_0(x)$. This convergence extends to multi-phase initial data in the distributional sense \cite{Jin-Xiantao}.

\subsection{One-Dimensional Moment System for Multi-Phase Solutions}
\label{Sec1}

For clarity of presentation, we focus on the one-dimensional case in physical and velocity space. The multi-phase solution ansatz for the Vlasov equation \eqref{Vlasov-eqn} takes the form:
\begin{equation}\label{ansatz} 
w(t,x,v) = \sum_{k=1}^{N} \rho_k(t,x) \delta(v-u_k(t,x)),
\end{equation}
where $N$ represents the number of distinct phases. The corresponding moments are defined as:
\begin{equation}\label{m-def} 
m_l(t,x) = \int_{\mathbb{R}} w(t,x,v) v^l \, dv, \quad l = 0,\ldots,2N,
\end{equation}
with the macroscopic density and velocity given by:
\begin{equation}
\rho(t,x) = m_0, \quad u(t,x) = \frac{m_1}{m_0}.
\end{equation}
Applying the moment method to \eqref{Vlasov-eqn} yields the following system of $2N$ equations:
\begin{equation}\label{moment_system}
\begin{cases}
\partial_t m_0 + \partial_x m_1 = 0, \\[2pt]
\partial_t m_1 + \partial_x m_2 = -m_0 \partial_x\Phi, \\[2pt]
\quad \vdots \\[2pt]
\partial_t m_{2N-1} + \partial_x m_{2N} = -(2N-1)m_{2N-2} \partial_x\Phi.
\end{cases}
\end{equation}
The ansatz \eqref{ansatz} establishes an explicit relationship between moments and phase variables:
\begin{equation}\label{m-condition}
m_l = \sum_{k=1}^N \rho_k u_k^l, \quad l = 0,\ldots,2N.
\end{equation}
For known $N$, the system \eqref{m-condition} can be inverted to express the highest-order moment as:
\begin{equation}\label{F-relation}
m_{2N} = F(m_0,\ldots,m_{2N-1}),
\end{equation}
completing the closure of \eqref{moment_system}. This inversion is guaranteed by the weak hyperbolicity and invertibility properties established in \cite{Jin-Xiantao}, which permit the unique determination of $(\rho_k, u_k)$ from the lower-order moments.

\subsection{Motivation and Goal of This Work}

When the number of phases $N$ is small ($N<5$) and the solution ansatz is known (as discussed in Section~\ref{Sec1}), one can explicitly relate $m_{2N}$ to lower-order moments $\{m_l\}_{l=0}^{2N-1}$ through the formula $m_l = \sum_{k=1}^N \rho_k u_k^l$. However, practical applications present two significant challenges:
\begin{enumerate}
    \item the value of $N$ is typically unknown, and
    \item for large $N$, inverting the resulting $2N\times 2N$ system to determine the closure relation $F$ in \eqref{F-relation} becomes computationally intractable.
\end{enumerate} These limitations motivate our deep learning approach to moment closure, enabling robust multi-phase computations for the semiclassical Schr\"odinger equation.

To illustrate the core challenge, consider the simplest case when $N=1$. The corresponding $2\times 2$ moment system takes the form:
\begin{equation}
\label{2-2moment}
\begin{cases}
\partial_t m_0 + \partial_x m_1 = 0, \\[4pt] 
\partial_t m_1 + \partial_x \mathcal{F} = -m_0\partial_x\Phi,
\end{cases}
\end{equation}
where $\mathcal{F}$ in the flux term no longer corresponds to the standard second-order moment from \eqref{m-def}. Our primary objective is to \textit{learn} the closure relation $\mathcal{F}$ that: (i) properly closes the system, and (ii) yields accurate approximations of $m_0$ and $m_1$ when compared to reference solutions of the Vlasov equation.

\subsection{Particle Method for the Vlasov Equation}

We employ a particle method \cite{FilbetPIC} to solve the Vlasov equation, representing the plasma distribution through a discrete set of computational particles. The method evolves particle trajectories along characteristic curves of the Vlasov equation:
\begin{equation}
\label{eq:characteristics}
\begin{cases}
\frac{dX}{dt} = V, \\[2pt]
\frac{dV}{dt} = -\nabla_x\Phi, \\[2pt]
X(0) = x^0, \quad V(0) = v^0.
\end{cases}
\end{equation}
The initial distribution $f_0$ is approximated by a sum of Dirac masses:
\begin{equation}
f_{M}^0(x,v) = \sum_{k=1}^{M} w_k \delta(x-x_k^0) \delta(v-v_k^0),
\end{equation}
where $\{(x_k^0, v_k^0)\}_{k=1}^M$ represents $M$ particles sampling the initial phase space density. The time-evolved solution then becomes:
\begin{equation}
f_M(t,x,v) = \sum_{k=1}^M w_k \delta(x-X_k(t)) \delta(v-V_k(t)),
\end{equation}
with $(X_k(t), V_k(t))$ following the characteristic equations \eqref{eq:characteristics}.

For numerical implementation, we regularize the Dirac masses using smooth shape functions:
\begin{equation}
f_{M,\alpha}(t,x,v) = \sum_{k=1}^M w_k \varphi_{\alpha}(x-X_k(t)) \varphi_{\alpha}(v-V_k(t)),
\end{equation}
where $\varphi_{\alpha}(z) = \alpha^{-d}\varphi(z/\alpha)$ is a scaled kernel satisfying:
\begin{itemize}
\item $\varphi \in C^r$ (r-times differentiable)
\item $\text{supp}(\varphi) \subset B_R(0)$ (compact support)
\item $\varphi(z) \geq 0$ (non-negative)
\item $\int_{\mathbb{R}^d} \varphi(z) dz = 1$ (normalized)
\item $\varphi(-z) = \varphi(z)$ (symmetric)
\end{itemize}
Common choices include B-splines and Gaussian kernels \cite{Splines}, with $\alpha$ controlling the smoothing length scale.

\section{Neural Network Approaches}\label{sec.2}

\subsection{A Two-Stage Moment Closure Model}
In this subsection, we will introduce our novel two-stage method to close the moment system of the Liouville equation, known as the semiclassical limit of the Schr\"odinger equation. 

\subsubsection{Physics Informed Neural Networks (PINNs)}
We begin by briefly reviewing Physics-Informed Neural Networks (PINNs) \cite{PINN1}. PINNs approximate solutions to partial differential equations (PDEs) by representing them with fully connected neural networks. A fully connected neural network with $L$ layers can be expressed recursively as follows for the $l$-th and $(l+1)$-th layers ($l = 1, 2, \dots, L-1$):
\begin{equation}\label{two-layer}
    n_j^{(l+1)} = \sum_{i=1}^{m_l} w_{ji}^{(l+1)} \sigma_l(n_i^{(l)}) + b_j^{(l+1)},
\end{equation}
where $n_i^{(l)}$ denotes the $i$-th neuron in the $l$-th layer, $\sigma_l$ is the activation function, $w_{ji}^{(l+1)}$ and $b_j^{(l+1)}$ represent the weights and biases for the $(l+1)$-th layer, and $m_l$ is the number of neurons in the $l$-th layer. The input layer $\{n_i^{(0)}\}_{i=1}^{m_0}$ typically represents physical variables, e.g., $(t, \boldsymbol{x})$.

Consider the general form of a PDE system:
\begin{equation}\label{eq:gen_pde}
    \begin{cases}
    \mathcal{D}[u(t,\boldsymbol{x}), t, \boldsymbol{x}] = 0, \quad (t, \boldsymbol{x}) \in [0, T] \times \Omega, \\
    \mathcal{I}[u(t,\boldsymbol{x}), t, \boldsymbol{x}] = 0, \quad (t, \boldsymbol{x}) \in \{0\} \times \Omega, \\
    \mathcal{B}[u(t,\boldsymbol{x}), t, \boldsymbol{x}] = 0, \quad (t, \boldsymbol{x}) \in [0, T] \times \partial\Omega,
    \end{cases}
\end{equation}
where $\mathcal{D}$ is the differential operator of the PDE, and $\mathcal{I}$ and $\mathcal{B}$ denote the initial and boundary condition operators, respectively. Here, $\Omega \subseteq \mathbb{R}^d$ is the spatial domain. PINNs approximate the solution $u(t, \boldsymbol{x})$ with a neural network $\hat{u}(t, \boldsymbol{x}; \theta)$. The key idea is to incorporate the residuals from the governing equations and physical constraints into the loss function:
\begin{multline}
\label{eq:pinn_prob_gen}
    \min_{\theta \in \mathbb{R}^p} L(\theta) := \frac{1}{N} \sum_{j=1}^N \left( \mathcal{D}[\hat{u}(x_r^j; \theta), x_r^j] \right)^2 + \frac{1}{M} \sum_{k=1}^M \left( \mathcal{I}[\hat{u}(x_i^k; \theta), x_i^k] \right)^2 \\+ \frac{1}{P} \sum_{l=1}^P \left( \mathcal{B}[\hat{u}(x_b^l; \theta), x_b^l] \right)^2,
\end{multline}
where $x_r^j$, $x_i^k$, and $x_b^l$ are the selected collocation points in different domains, while $N$, $M$, and $P$ are the number of the collocation points.  Here, automatic differentiation is used to compute the PDE residuals at a set of collocation points, enforcing physical consistency during training. The parameters $\theta$ are optimized via stochastic gradient descent or adaptive methods such as Adam \cite{Adam14}. After training, the resulting surrogate $\hat{u}(t, \boldsymbol{x}; \theta^\star)$ can be evaluated at any point in the domain.

As a motivating example, consider the single-species case ($N=1$) in the system \eqref{2-2moment}. Our goal is to approximate the flux term $\partial_x \mathcal{F}$ using training data. Details on how this term is constructed will be discussed later. Assuming an $L$-layer neural network, the input is $(t, \boldsymbol{x})$, and the output of the final layer is used to represent $m_0^{\text{NN}}(t, \boldsymbol{x}; \theta)$ and $m_1^{\text{NN}}(t, \boldsymbol{x}; \theta)$, where $\theta$ denotes the learnable parameters.
 
\subsubsection{Model Design}
After trials and errors, we find that the following architecture is the most efficient. We first learn the relation between the spatial derivative of the last moment with all the previous moments, then adopt PINNs to solve the system \eqref{moment_system}. For simplicity, we focus on the case $N=1$, which can be directly extended to larger $N$ (multi-phase case). We summarize the framework of our method for the forward problem in Figure \ref{fig:moment_framework} below, here $\sigma(x)$ denotes the activation function such as the hyperbolic tangent function. 
\begin{figure}[htbp]
\centerline{\includegraphics[width=1\linewidth]{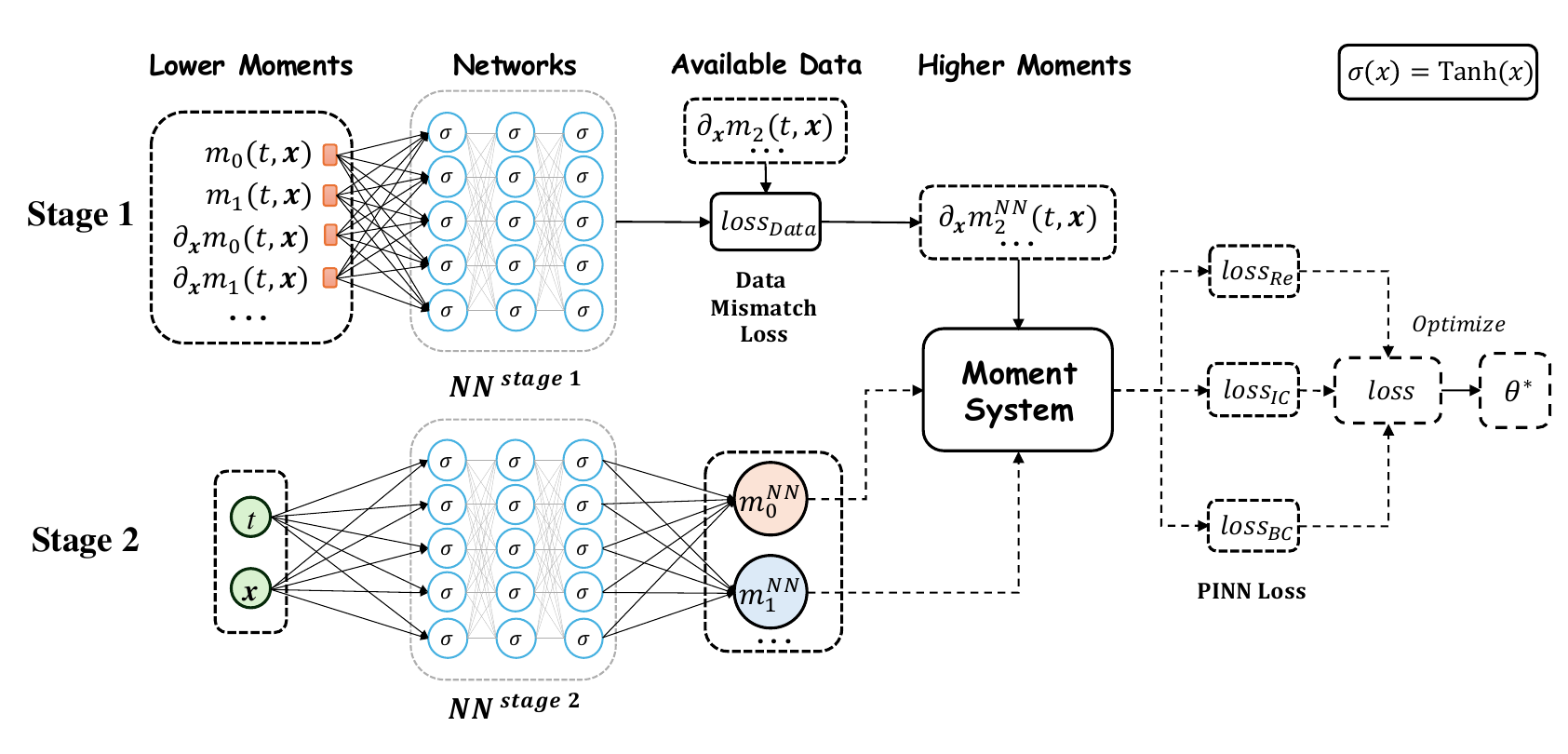}}
    \caption{Framework of the two-stage moment closure model of solving moment system.}
    \label{fig:moment_framework}
\end{figure}

\subsubsection{Loss Functions}
Considering the one-dimensional spatial domain $x\in\Omega$ and time $t\in[0,T]$ scenario, we would like to introduce the proposed method for both stages, including the design of the loss function.

\noindent \textbf{Stage 1.} Let the output of neural network (NN) be $\partial_x m_2^{NN}$, and inputs of NN are $m_0$, $m_1$, $\partial_x m_0$, $\partial_x m_1$. In Stage 1, we study the relation between $m_0$, $m_1$, $\partial_x m_0$, $\partial_x m_1$ and $\partial_x m_2$ by using the neural network: 
$$ \partial_x m_2^{\text{NN},\textit{stage1}} := \mathcal{NN}(m_0, m_1, \partial_x m_0, \partial_x m_1). $$
The loss in this stage is defined by
\begin{equation}
\label{S1}
Loss_{\textit{stage1}} := \|\partial_x m_2^{NN} - \partial_x m_2^{Data}\|_{L^2_{t,x}} = \sum_{i=1}^{K} || \partial_x m_2^{NN}(t_i, x_i) - \partial_x m_2^{Data}(t_i, x_i) ||^2. 
\end{equation}

\noindent \textbf{Stage 2.} 
Design two neural networks where the inputs are $t$, $x$, and outputs are the approximations for $m_0(t,x)$ and $m_1(t,x)$. We employ the PINNs on the moment system \eqref{2-2moment}:   \begin{equation}
    \label{2-2moment_sp2}
      \begin{cases}
      \partial_t m_0 + \partial_x m_1 = 0,  \\[6pt]
      \partial_t m_1 + \partial_x m_2^{\text{NN},\textit{stage1}}  = -m_0 \partial_x \Phi,  \end{cases}
\end{equation} 
where $\partial_x m_2^{\text{NN},\textit{stage1}}$ is obtained from Stage 1. To find the optimal values for the network parameters $\theta$ that are composed of all the weights and biases, the neural network is trained by minimizing the loss function 
\begin{equation}
\label{Loss-PINN}
\begin{aligned}
Loss_{\textit{stage2}}(\theta) = Loss_{\text{GE}}(\theta) + \lambda_1 Loss_{\text{BC}}(\theta) + \lambda_2 Loss_{\text{IC}}(\theta),  
 \end{aligned}
\end{equation}
where the loss for the residual of governing equations 
\eqref{2-2moment_sp2} is composed of two parts: 
\begin{align}
& Loss_{\text{GE}_1}(\theta) = \int_{\mathcal{T}}\int_{\Omega}\|\partial_t m_0^{\text{NN}} + \partial_x m_1^{\text{NN}}\|^2 dx dt. \label{L-GE1} \\
& Loss_{\text{GE}_2}(\theta) = \int_{\mathcal{T}}\int_{\Omega} \|
\partial_t m_1^\text{NN} + \partial_x m_2^{\text{NN},\textit{stage1}}  + m_0^\text{NN} \partial_x \Phi\|^2 dx dt, \label{L-GE2}
\end{align}
and $Loss_{\text{GE}}= Loss_{\text{GE}_1} + Loss_{\text{GE}_2}$. 
Assume periodic boundary condition (BC), the loss induced by BC is given by 
\small{\begin{equation}\label{L-BC}
\begin{aligned}
Loss_{\text{BC}}(\theta) = & \int_{\mathcal{T}}\|m_0^{\text{NN}}(s,x=x_L) - m_0^{\text{NN}}(s,x=x_R)\|^2 ds \\ 
& + \int_{\mathcal{T}}\|m_1^{\text{NN}}(s,x=x_L) -  m_1^{\text{NN}}(s,x=x_R)\|^2 ds, 
\end{aligned}
\end{equation}}
and the loss of the initial condition is  
\small{\begin{equation}\label{L-IC}
\begin{aligned}
Loss_{\text{IC}}(\theta) = & \int_{\Omega} \| m_0^{\text{NN}}(t=0,\tau)-m_0(t=0,\tau) \|^2 d\tau  \\
& + \int_{\Omega} \| m_1^{\text{NN}}(t=0,\tau)-m_1(t=0,\tau) \|^2 d\tau. 
\end{aligned}
\end{equation}}

In order to determine whether our designed neural network-based moment closure system is efficient and accurate, we compute the approximated solutions $m_0^{NN}$, $m_1^{NN}$ and compare it with the reference solution computed by PIC or finite-volume method for the Vlasov equation. 

\subsection{Convergence of the Loss and Solutions for Stage 2}

In Stage 2 and the framework of PINNs, we adapt the Universal Approximation Theorem (UAT) \cite[Theorem 2.1]{li1996simultaneous} to our moment system \eqref{2-2moment} and show the convergence of loss function and NN approximated solutions, where $\partial_x m_2^{\text{NN},\textit{stage1}}$ is given from Stage 1.  
\begin{lemma}[Li, Theorem 2.1, \cite{li1996simultaneous}]
\label{Lemma1}
Suppose a solution $(m_0,m_1)$ of the moment system 
\eqref{2-2moment} satisfies $m_0,\ m_1 \in C^1([0,T]\times \Omega)$. Let the activation function $\sigma$ be any non-polynomial function in $C^1(\mathbb R)$. Then for any $\delta>0$ there exists a two-layer neural network for $k=0$ and $1$,
$$ m_k^{\text{NN}}(t,x) = \sum_{i=1}^{n_{k,1}} w_{1i}^{k,(2)} \sigma\left( \left(w_{i 1}^{k,(1)}, w_{i 2}^{k,(1)}\right) \cdot(t, x)+b_{i}^{k,(1)}\right) + b_{i}^{k,(2)}, 
$$
such that 
\begin{equation}
\label{UAT}
\begin{aligned}
\displaystyle &  \left\| m_k - m_k^{\text{NN}} \right\|_{L^{\infty}(K)} < \delta, \qquad 
\left\| \partial_t ( m_k - m_k^{\text{NN}}) \right\|_{L^{\infty}(K)} < \delta, \\[4pt]
& \left\| \partial_x ( m_k - m_k^{\text{NN}}) \right\|_{L^{\infty}(K)} < \delta, 
\end{aligned}
\end{equation}
where the domain $K$ denotes $[0,T]\times\Omega$. 
\end{lemma}
The above result can be generalized to neural network with several hidden layers \cite{Multilayer}. We begin by stating our first main theorem, which guarantees that if the equation admits a \(C^1\) solution, then there exists a sequence of neural network solutions whose total loss term converges to zero:
\begin{theorem}\label{forward_loss}Denote $K=[0,T]\times\Omega$. Assume that the solution $(m_0,m_1)$ to the moment system \eqref{2-2moment} with the initial profiles $(m_0(0,\cdot),m_1(0,\cdot))$ and the periodic boundary conditions belongs to $C^1(K)$, and the activation function $\sigma\in C^1(K)$ is non-polynomial. Then, for both k=0 and 1, there exist neural network weights $\{n_{k,[j]}, w_{k,[j]}, b_{k,[j]}\}_{j=1}^\infty$ such that a sequence of the DNN solutions with $n_{k,[j]}$ nodes, denoted by $$\{(m_k)_j(t,x) = m_k^{\text{NN}}(t,x;n_{k,[j]}, w_{k,[j]}, b_{k,[j]})\}_{j=1}^{\infty},$$ satisfies
\begin{equation}\label{forward_loss_0}
\text{Loss}_{stage2}(\theta)\rightarrow 0\text{ as }j\rightarrow\infty.
\end{equation}
\end{theorem} 
\begin{proof} Fix any $\delta>0$. By Lemma \ref{Lemma1}, we have the existence of neural network solutions $(m_k)_j$ such that \eqref{UAT} holds for both $k=0$ and $k=1$. Then note that by \eqref{UAT} and by the fact that $m_0$ and $m_1$ are solutions to the system \eqref{2-2moment_sp2}, we have
\small{\begin{equation*}
\begin{aligned}
   & Loss_{\text{GE}}(\theta)\\
    &=Loss_{\text{GE}_1}(\theta)+Loss_{\text{GE}_2}(\theta)\\[2pt]
   & =\int_K\left(|\partial_t (m_0)_j + \partial_x (m_1)_j|^2+|\partial_t (m_1)_j + \partial_x m_2^{\text{NN},\textit{stage1}}  +(m_0)_j \partial_x \Phi|^2\right)dxdt\\[2pt]
   & =\int_K\left(|\partial_t (m_0-(m_0)_j) + \partial_x (m_1-(m_1)_j)|^2+|\partial_t (m_1-(m_1)_j) +(m_0-(m_0)_j )\partial_x \Phi|^2\right)dxdt\\[2pt]
   & \le (4\delta^2 +(\delta+\delta \|\partial_x \Phi\|_{L^\infty_x(\Omega)})^2)T\mu(\Omega),
    \end{aligned}
    \end{equation*}}
where $\mu$ is the Lebesgue measure on $\mathbb{R}.$ In addition, the loss for the boundary data is bounded from above by 
\small{\begin{equation*}
\begin{aligned}
Loss_{\text{BC}}(\theta) 
&= \int_0^T\|(m_0)_j(t,x=x_L) - (m_0)_j(t,x=x_R)\|^2 dt  \\[2pt]
&\quad + \int_0^T\|(m_1)_j(t,x=x_L) -  (m_1)_j(t,x=x_R)\|^2 dt\\[2pt]
&=\int_0^T\|(m_0-(m_0)_j)(t,x=x_L) - (m_0-(m_0)_j)(t,x=x_R)\|^2 dt \\[2pt]
&\quad + \int_0^T\|(m_1-(m_1)_j)(t,x=x_L) -  (m_1-(m_1)_j)(t,x=x_R)\|^2 dt \le 8\delta^2 T,
 \end{aligned}
    \end{equation*}}
since $m_0$ and $m_1$ are periodic in $x.$ Moreover, the loss for the initial condition goes by 
\begin{equation*}
\begin{aligned}
Loss_{\text{IC}}(\theta) & = \int_{\Omega} \| (m_0)_j(t=0,x)-m_0(t=0,x) \|^2 dx  \\
 &\quad + \int_{\Omega} \| (m_1)_j(t=0,x)-m_1(t=0,x) \|^2 dx \le 4\delta^2\mu(\Omega).
 \end{aligned}
    \end{equation*}
    Therefore, the total loss $$\text{Loss}_{\textit{stage2}}(\theta) = Loss_{\text{GE}}(\theta) + \lambda_1 Loss_{\text{BC}}(\theta) + \lambda_2 Loss_{\text{IC}}(\theta)$$ can be made arbitrarily small by choosing sufficiently small $\delta>0.$ Since for each $\delta=\frac{1}{j}>0$ we can construct the corresponding $j^{th}$ sequential elements $(m_0)_j$ and $(m_1)_j$ such that \eqref{UAT} holds with $\delta =\frac{1}{j}. $ As $j\to +\infty,$ we have $\delta \to 0^+$ and obtain that the total loss $\text{Loss}_{\textit{stage2}}(\theta) $ vanishes. This completes the proof.
\end{proof}

We note that Theorem \ref{forward_loss} ensures the existence of neural network weights that can reduce the error function to an arbitrary level. However, the convergence of the loss function to zero does not necessarily imply that the neural network solution converges to the true solution of the original equation. 

To address this limitation, we present our second main result, Theorem \ref{theorem_forward}, which establishes that the neural network architecture converges to a $C^1$ solution in an appropriate function space when weights minimize \( Loss_{Total} \). 
\begin{theorem}\label{theorem_forward} Let $n,$ $w$, and $b$ be the number of nodes, the weights, and the biases for the neural network architecture. Suppose that $\partial_x m_2^{\text{NN},\textit{stage1}} $ is obtained from Stage 1. Denote the DNN solutions $m_k^{\text{NN}}$ for each $k=0$ and $k=1$ as $m_k^{\text{NN}}(t,x) = m_k^{\text{NN}}(t,x;n, w, b)$.  Define the energy functional as
\[
E(t) = \frac{1}{2} \int_{\Omega} \big((m_0^{\text{NN}}(t,x)-m_0(t,x))^2 + (m_1^{\text{NN}}(t,x)-m_1(t,x))^2\big) \, dx,
\]
where $(m_0,m_1)$ is the solution to  \eqref{2-2moment} with initial profiles $(m_0(0,\cdot),m_1(0,\cdot))$ under the periodic boundary condition. Then, for any $t\in[0,T],$ $E(t)\to 0$ if $$\text{Loss}_{stage2}(\theta)\rightarrow 0 \text{ and }
\|m_0^{\text{NN}}-m_0\|_{L^1([0,t]; \dot{H}^1(\Omega))} \to 0.$$
\end{theorem}
\begin{proof}
Define 
\begin{equation}
\begin{aligned}
& D_1(t,x)=\partial_t m_0^{\text{NN}} + \partial_x m_1^{\text{NN}},\\[4pt]
& D_2(t,x) = 
\partial_t m_1^{\text{NN}} + \partial_x m_2^{\text{NN},\textit{stage1}}  + C_1(x)m_0^{\text{NN}} , \\[4pt]
& D_3(t) =  |m_0^{\text{NN}}(t,x=x_L) - m_0^{\text{NN}}(t,x=x_R)|^2 +
|m_1^{\text{NN}}(t,x=x_L) -  m_1^{\text{NN}}(t,x=x_R)|^2,  \\[4pt]
 & D_4(x) =  |m_0^{\text{NN}}(t=0,x)-m_0(t=0,x)|^2+| m_1^{\text{NN}}(t=0,x)-m_1(t=0,x) |^2,
\end{aligned}
\end{equation}
where we denote $C_1(x)=\partial_x \Phi(x).$
Since  $(m_0,m_1)$ be a solution to  \eqref{2-2moment} with the initial profiles $(m_0(0,\cdot),m_1(0,\cdot))$ under the periodic boundary condition, we obtain that their difference from the neural network solutions $(m_0^{\text{NN}},m_1^{\text{NN}})$ solve the following system of equations:\begin{equation}\label{system of equations for differences}
\begin{cases}
\partial_t (m_0^{\text{NN}} - m_0) + \partial_x (m_1^{\text{NN}} - m_1) = D_1(t,x), \\
\partial_t (m_1^{\text{NN}} - m_1) + C_1(x) (m_0^{\text{NN}} - m_0) = D_2(t,x).
\end{cases}
\end{equation}We assume that \(D_1, D_2 \in L^2_{t,x}\) and both converge to zero in \(L^2_{t,x}\) as $j\to 0$ by the previous theorem (Theorem \ref{forward_loss}). Now for each $t\in [0,T]$, define the energy functional as
\[
E(t) = \frac{1}{2} \int_{\Omega} \big((m_0^{\text{NN}}(t,x)-m_0(t,x))^2 + (m_1^{\text{NN}}(t,x)-m_1(t,x))^2\big) \, dx,
\]
Formally, differentiating with respect to time, we observe
\[
\frac{d}{dt} E(t) = \int_{\Omega} \big((m_0^{\text{NN}}-m_0)\partial_t (m_0^{\text{NN}}-m_0) + (m_1^{\text{NN}}-m_1)\partial_t (m_1^{\text{NN}}-m_1)\big) \, dx.
\]
Using the system \eqref{system of equations for differences} and substituting the derivatives into the energy derivative, we obtain
\begin{equation*}
\begin{aligned}
    \frac{d}{dt} E(t) & = \int_{\Omega} \big[(m_0^{\text{NN}}-m_0)(-\partial_x (m_1^{\text{NN}}-m_1) + D_1) \\
    &\quad + (m_1^{\text{NN}}-m_1)(-C_1(x)(m_0^{\text{NN}}-m_0) + D_2)\big] dx.
\end{aligned}
\end{equation*}

Now we simplify the energy derivative by handling each term. Firstly, by taking the integration by parts on the first term, we obtain
\[
\int_{\Omega} (m_0^{\text{NN}}-m_0)(-\partial_x (m_1^{\text{NN}}-m_1)) \, dx = \int_{\Omega} (\partial_x (m_0^{\text{NN}}-m_0)) (m_1^{\text{NN}}-m_1) \, dx+D_{\text{boundary}},
\]where by the Cauchy-Schwarz inequality, we have 
$|D_{\text{boundary}}|\le D_3(t)$. 
Therefore, 
\begin{equation*}
\begin{aligned}
    \frac{d}{dt} E(t) & \le D_3(t)+ \int_{\Omega} (\partial_x (m_0^{\text{NN}}-m_0))(m_1^{\text{NN}}-m_1) \, dx + \int_{\Omega} (m_0^{\text{NN}}-m_0) D_1 \, dx\\[2pt]
    &\quad + \int_{\Omega} (m_1^{\text{NN}}-m_1) D_2 \, dx - \int_{\Omega} C_1(x) (m_0^{\text{NN}}-m_0)(m_1^{\text{NN}}-m_1) \, dx.
\end{aligned}
\end{equation*}
By further using Cauchy–Schwarz and Young’s inequalities, we estimate each component as follows:
\begin{itemize}
    \item \(\displaystyle\left| \int (\partial_x (m_0^{\text{NN}}-m_0))(m_1^{\text{NN}}-m_1) \, dx \right| \leq \|\partial_x (m_0^{\text{NN}}-m_0)\|_{L^2_x} \|(m_1^{\text{NN}}-m_1)\|_{L^2_x}\).
    \item \(\displaystyle\left| \int (m_0^{\text{NN}}-m_0) D_1 \, dx \right| \leq \|(m_0^{\text{NN}}-m_0)\|_{L^2_x} \|D_1\|_{L^2_x}\).
    \item \(\displaystyle\left| \int (m_1^{\text{NN}}-m_1) D_2 \, dx \right| \leq \|(m_1^{\text{NN}}-m_1)\|_{L^2_x} \|D_2\|_{L^2_x}\).
    \item \(\displaystyle\left| \int C_1(x)(m_0^{\text{NN}}-m_0)(m_1^{\text{NN}}-m_1) \, dx \right| \leq \|C_1\|_{L^\infty} \|(m_0^{\text{NN}}-m_0)\|_{L^2_x} \|(m_1^{\text{NN}}-m_1)\|_{L^2_x}\).
\end{itemize}
Apply the Young’s inequality, 
we obtain
\[
\frac{d}{dt} E(t) \leq D_3(t)+ C (E(t)+\|\partial_x (m_0^{\text{NN}}-m_0)\|_{L^2_x}^2) + \frac{1}{2}\|D_1(t,\cdot)\|_{L^2_x}^2 + \frac{1}{2}\|D_2(t,\cdot)\|_{L^2_x}^2,
\]
for some constant \(C\) depending on \(\|C_1\|_{L^\infty}\).
Integrating it from \(0\) to \(s\) for any $s\in[0,T]$, we have
\begin{equation*}
\begin{aligned}
    E(s) &\leq E(0) +  \int_0^s (D_3(t)+C\|\partial_x (m_0^{\text{NN}}-m_0)\|_{L^2_x}^2+CE(t)) \, dt \\[2pt]
    &\quad + \frac{1}{2} \int_0^s \|D_1(t,\cdot)\|_{L^2_x}^2 \, dt + \frac{1}{2} \int_0^s \|D_2(t,\cdot)\|_{L^2_x}^2 \, dt, 
\end{aligned}
\end{equation*}
where we use that
$E(0)=\int_{\Omega}D_4(x)dx$. By applying Gr\"onwall’s inequality, one finally obtains that for any $s\in[0,T]$, 
\begin{equation*}
\begin{aligned}
    E(s) &\leq \bigg(\int_{\Omega}D_4(x)dx+\int_0^s (D_3(t)+C\|\partial_x (m_0^{\text{NN}}-m_0)\|_{L^2_x}^2) \, dt\\[2pt] &\qquad +\frac{1}{2}\|D_1\|_{L^2_{t,x}}^2 + \frac{1}{2}\|D_2\|_{L^2_{t,x}}^2 \bigg)\,  e^{C s}.
\end{aligned}
\end{equation*}
 By definition, if $\text{Loss}_{stage2}(\theta)\rightarrow 0$, then $$\left(\int_{\Omega}D_4(x)dx+\int_0^s D_3(t) \, dt+\frac{1}{2}\|D_1\|_{L^2_{t,x}}^2 + \frac{1}{2}\|D_2\|_{L^2_{t,x}}^2 \right)\to 0.$$ Therefore, we conclude that for each  $s\in[0,T]$, $E(s)\to 0$ as $$\text{Loss}_{stage2}(\theta)\rightarrow 0 \text{ and }
\|m_0^{\text{NN}}-m_0\|_{L^1([0,t]; \dot{H}^1(\Omega))} \to 0.$$ This completes the proof.
\end{proof}
\begin{remark}
    The ``pathological" term $\|\partial_x (m_0^{\text{NN}}-m_0)\|_{L^2_x}$ does not vanish and cannot be absorbed into the energy without further structure. To control it, one may need second-level energy involving $\partial_x \left( m_0^{\text{NN}}-m_0 \right)$ or some dissipation from higher-order terms (if the original system has them). Without such mechanisms, this term may cause solution growth, and energy decay alone is insufficient to ensure convergence.
\end{remark}

\section{Numerical Experiments}\label{sec.3}
In our numerical experiments, we will show several examples to illustrate the effectiveness of our designed \textit{two-stage moment closure model}. For the first two numerical examples, we employ our method for the moment equation \eqref{2-2moment}. For the third high-dimensional example, we employ the method for the moment equation \eqref{test3_system}. In the first stage, 
we evaluate the effectiveness of our novel model-free learning approach in acquiring $\partial_x m_2$ against some alternative learning methodologies \cite{huang2022machine}. Our method outperforms the others due to its network architecture, which has a more flexible structure and incorporates a broader range of information. During the second stage, we adapt the PINNs to obtain the approximation of lower moments such as $m_0$ and $m_1$ using learned moments in stage $1$.

\textbf{Learning Schemes.} We compare four different learning schemes in Stage 1. The training inputs for Stage 1 vary with each scheme. \textit{Neural augmented scheme}  is our preferred method, which takes the form of $\partial_x m_2^{NN} = \mathcal{NN}(m_0, m_1, \partial_x m_0, \partial_x m_1)$, with training sets $m_0$, $m_1$, $\partial_x m_0$, and $\partial_x m_1$ obtained by the traditional numerical method.
In \textit{derivative-only neural scheme}, we consider the form $\partial_x m_2^{NN} = \mathcal{NN}(\partial_x m_0, \partial_x m_1)$, with the training inputs being $\partial_x m_0$ and $\partial_x m_1$. \textit{Linear neural scheme} utilizes a linear combination form similar to \cite{huang2022machine}, where $\partial_x m_2^{NN} = {c_0}^{\mathcal{NN}}\partial_x m_0 + {c_1}^{\mathcal{NN}}\partial_x m_1$, with the training inputs being $\partial_x m_0$, $\partial_x m_1$, and ${c_0}^{\mathcal{NN}}$, ${c_1}^{\mathcal{NN}}$ being the coefficients learned by neural networks, where ${c_0}^{\mathcal{NN}} = \mathcal{NN}(m_0,m_1)$. In \textit{extended linear neural scheme}, we study a more complicated linear combination in the form of $\partial_x m_2^{NN} = {c_0}^{\mathcal{NN}}m_0 + {c_1}^{\mathcal{NN}}m_1 + {c_2}^{\mathcal{NN}} \partial_x m_0 + {c_3}^{\mathcal{NN}} \partial_x m_1$, with the training inputs being $m_0$, $m_1$, $\partial_x m_0$ and $\partial_x m_1$. Here, the coefficients ${c_0}^{\mathcal{NN}}$, ${c_1}^{\mathcal{NN}}$, ${c_2}^{\mathcal{NN}}$ and ${c_3}^{\mathcal{NN}}$ are learned by neural networks. The training inputs for Test I and Test II are different from Test III, which will be introduced in section \ref{Test3}. During the training in Stage 2, we only need the datasets $\partial_x m_2^{NN}$ learned in Stage 1 and there is no need for other training datasets.
For notation simplicity, we conclude the schemes with their abbreviations in Table \ref{tab:scheme_summary}.
\begin{table}[htbp]
    \centering
    \begin{small}
            \setlength{\tabcolsep}{5.5pt}
                \begin{tabular}{lccc}
                    \toprule
                    Scheme & Abbreviation \\
                    \midrule
                     Neural Augmented Scheme (\textbf{Ours, Preferred}) & \textbf{scheme 1} & \\
                    Derivative-Only Neural Scheme (Ours)  & scheme 2 \\
                    Linear Neural Scheme (Huang et al., \cite{huang2022machine}) & scheme 3 \\
                    Extended Linear Neural  Scheme (Ours) & scheme 4 \\
                    \bottomrule
                \end{tabular}
            \caption{Schemes Summary.}
        \label{tab:scheme_summary}
    \end{small}
\end{table}

\textbf{Networks Architecture.}
In both Stage 1 and Stage 2, we approximate the solutions by the feed-forward neural network (FNN) with one input layer, one output layer, and $4$ hidden layers with $128$ neurons in each layer, unless otherwise specified. The hyperbolic tangent function (Tanh) is chosen as our activation function.

\textbf{Training Settings.}
The neural networks are trained by Adam with Xavier initialization. We set epochs to be $200000$ for Stage 1 and $50000$ for Stage 2, the learning rate as $10^{-3}$, and use full batch for most tests in the numerical experiments unless otherwise specified. All the hyper-parameters are chosen by trial and error. 

\textbf{Empirical Loss Design.}
For Test I and Test II, the empirical risk for stage $1$ is as follows:
\begin{equation}
\label{empirical_Loss-stage}
\begin{aligned}
\mathcal{R}_{\mathrm{Stage 1}} & 
 = \sum_{i=1}^{N_1} \left| \partial_x m_2^{\text{NN}}(t_i, x_i) - \partial_x m_2^{Data}(t_i, x_i) \right|^2.
\end{aligned}
\end{equation}
The empirical risk for stage $2$ is as follows:
\small{\begin{equation}
\label{empirical_Loss-PINN}
\begin{aligned}
\mathcal{R}_{\mathrm{Stage 2}} & = \frac{1}{N_1} \sum_{i=1}^{N_1} \Big|\partial_t m_0^{\text{NN}}(t_i,x_i) + \partial_x m_1^{\text{NN}}(t_i,x_i)\Big|^2 \,\\[4pt]
& + \frac{1}{N_2} \sum_{i=1}^{N_2} \Big| \partial_t m_1^\text{NN}(t_i,x_i) + \partial_x m_2^{\text{NN},\textit{stage1}}(t_i,x_i)  + m_0^\text{NN}(t_i,x_i) \partial_x \Phi(t_i,x_i)\Big|^2  \\[4pt]
& + \frac{\lambda_1}{N_3} \sum_{i=1}^{N_3} \left|\mathcal{B} \left( m_0^{\text{NN}}(t_i,x_i)  - m_0(t_i,x_i) \right)\right|^2 + \frac{\lambda_2}{N_3} \sum_{i=1}^{N_3} \left|\mathcal{B} \left( m_1^{\text{NN}}(t_i,x_i) - m_1(t_i,x_i)\right)\right|^2 \\[4pt]
& + \frac{\lambda_3}{N_4} \sum_{i=1}^{N_4} 
\left| m_0^{\text{NN}}(0,x_i) - m_0(0,x_i)\right|^2 + 
\frac{\lambda_4}{N_4}\sum_{i=1}^{N_4}\left|m_1^{\text{NN}}(0,x_i)) - m_1(0,x_i)\right|^2,
\end{aligned}
\end{equation}}
where $N_1$, $N_2$, $N_3$ $N_4$ are the number of sample points of $\mathcal{T} \times \Omega$, $\mathcal{T} \times \Omega$, $\mathcal{T} \times \partial \Omega$ and $\Omega$. 
For spatial points $x_i$, we select interior points evenly on $[x_l, x_r]$. For temporal points $t_i$, interior points are evenly picked in $[0, T]$. The tensor-product grids for the collocation points are used, and we set the penalty parameters in \eqref{empirical_Loss-PINN} as $(\lambda_1, \lambda_2,\lambda_3, \lambda_4) = (1,1,1,1)$.  For Test III, the empirical loss functions for stages $1$ and $2$ are shown in the Appendix. 

\textbf{Evaluation Matrices.} The reference solutions for moment quantities are obtained by the finite volume method and PIC method in solving the Vlasov equations. For Test I and Test II, we use the PIC method with the number of particles be $240000\times 300$. For Test III, we utilize the finite volume method with $\Delta t = 0.005$, $\Delta x = 0.02$, and $\Delta v = 0.2$. 
We compute relative $\ell^2$ 
errors for the moments between reference solutions and that approximated by our neural networks, with the relative $\ell^2$ for each time step error at time $t$ defined by: 
\begin{equation*}
\mathcal{E}_{\ell^2}:=\sqrt{\frac{\sum_{j=1}^N|m_{i}^{\mathcal{NN}}(t,\bx_j)-m_{i}^{\text{ref}}(t,\bx_j)|^2}{\sum_{j=1}^N|m_{i}^{\text{ref}}(t,\bx_j)|^2}}, \quad i=0,1,2. 
\end{equation*}
We also investigate mean square errors in some tests since using multiple matrices can reflect more comprehensive results. The mean square error at time $t$ is defined as follows:
\begin{equation*}
\mathcal{E}_{\text{MSE}}:= \frac{1}{N} \sum_{j=1}^N  |m_{i}^{\mathcal{NN}}(t, \bm{x}_j) - m_{i}^{\text{ref}}(t, \bm{x}_j)|^2, \quad i=0,1,2.
\end{equation*}
For Test I and Test II, we consider a one-dimensional spatial variable, for Test III, a two-dimensional problem is studied. Here $N$ is the number of total spatial points. 
Our experiments are conducted on a server with Intel$^{\text{(R)}}$ Xeon$^{\text{(R)}}$ Gold $6252$ and two A40 48GB GPUs.


\subsection{Test I: Moment closure for 1D Problem} 
Assume the initial data of the Vlasov equation \eqref{Vlasov-eqn} as
\begin{equation*} w(t=0,x,v) = \rho_0(x) \delta(v - \partial_x S_0), 
\end{equation*}
where 
\begin{equation*} \rho_0(x) = e^{-100(x-1)^2}, \quad 
S_0(x) = -\frac{1}{5}\ln\Big( e^{5(x-1)} + e^{-5(x-1)}\Big). 
\end{equation*}
The potential function is given by $\Phi=x^2/2$ (harmonic oscillator). Let the spatial domain be $[0,2]$ and velocity domain $[-2,2]$. We set $N_x=300$, $\Delta v=0.002$, $\Delta t=0.01$ and the final computational time $T=0.5$. Periodic boundary conditions are considered. 

\vspace{2mm}

\begin{figure}[htbp]
    \centering
    \begin{minipage}[b]{0.32\textwidth}
        \centering
        \includegraphics[width=\textwidth]{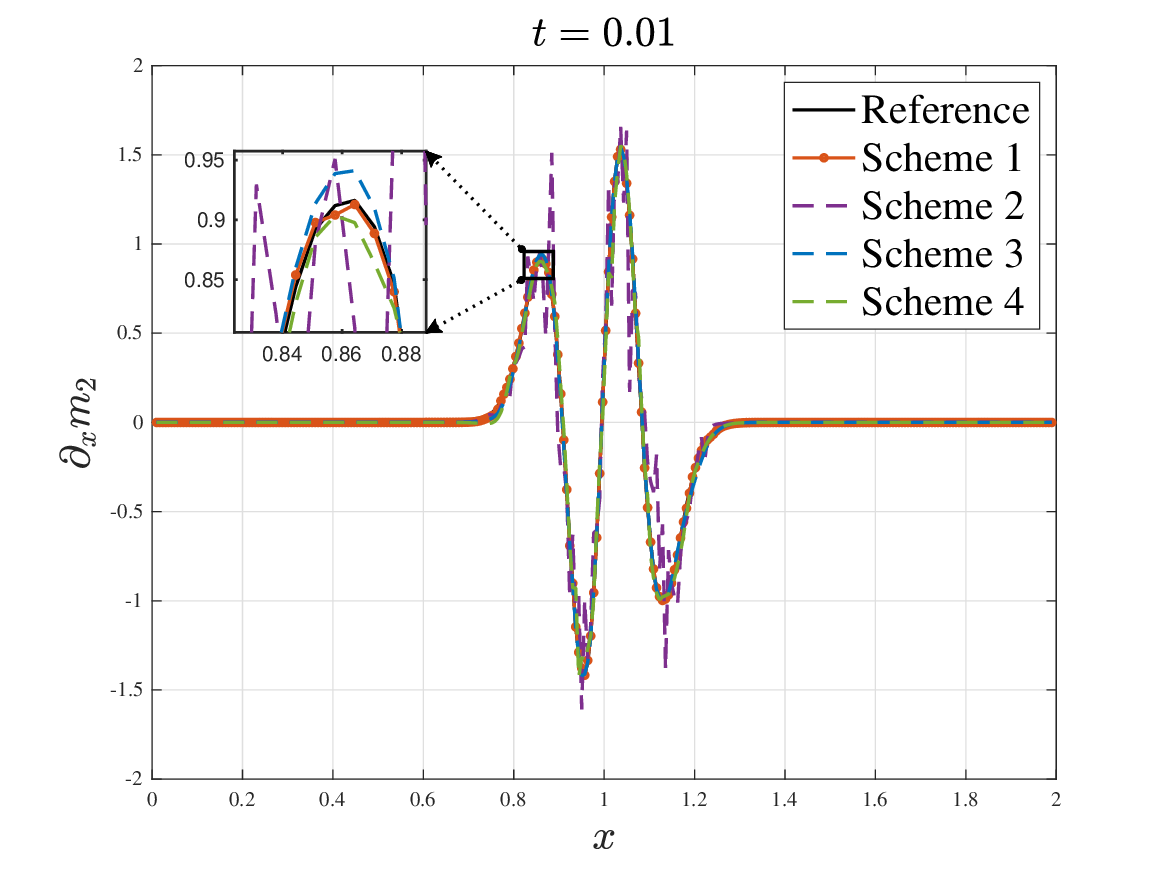} (a)
    \end{minipage}
    \begin{minipage}[b]{0.32\textwidth}
        \centering
        \includegraphics[width=\textwidth]{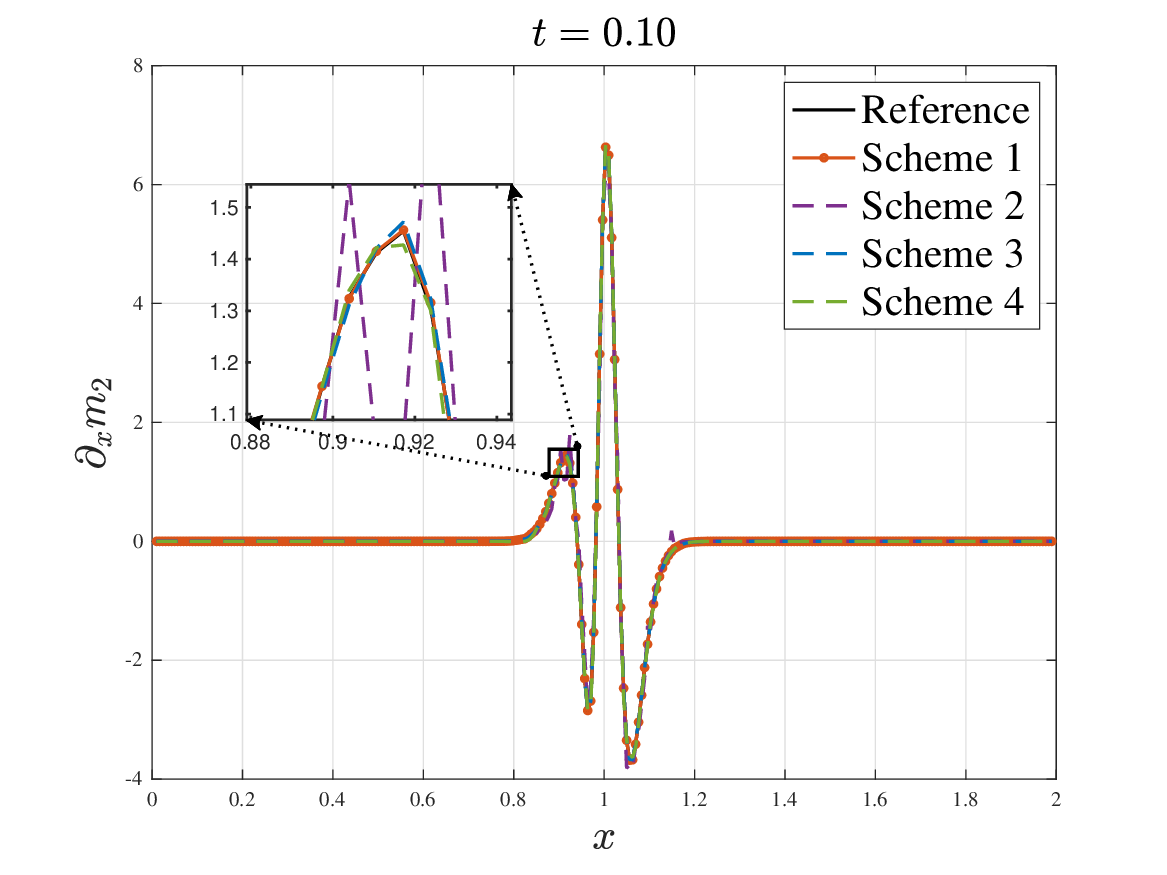} (b)
    \end{minipage}
    \begin{minipage}[b]{0.32\textwidth}
        \centering
        \includegraphics[width=\textwidth]{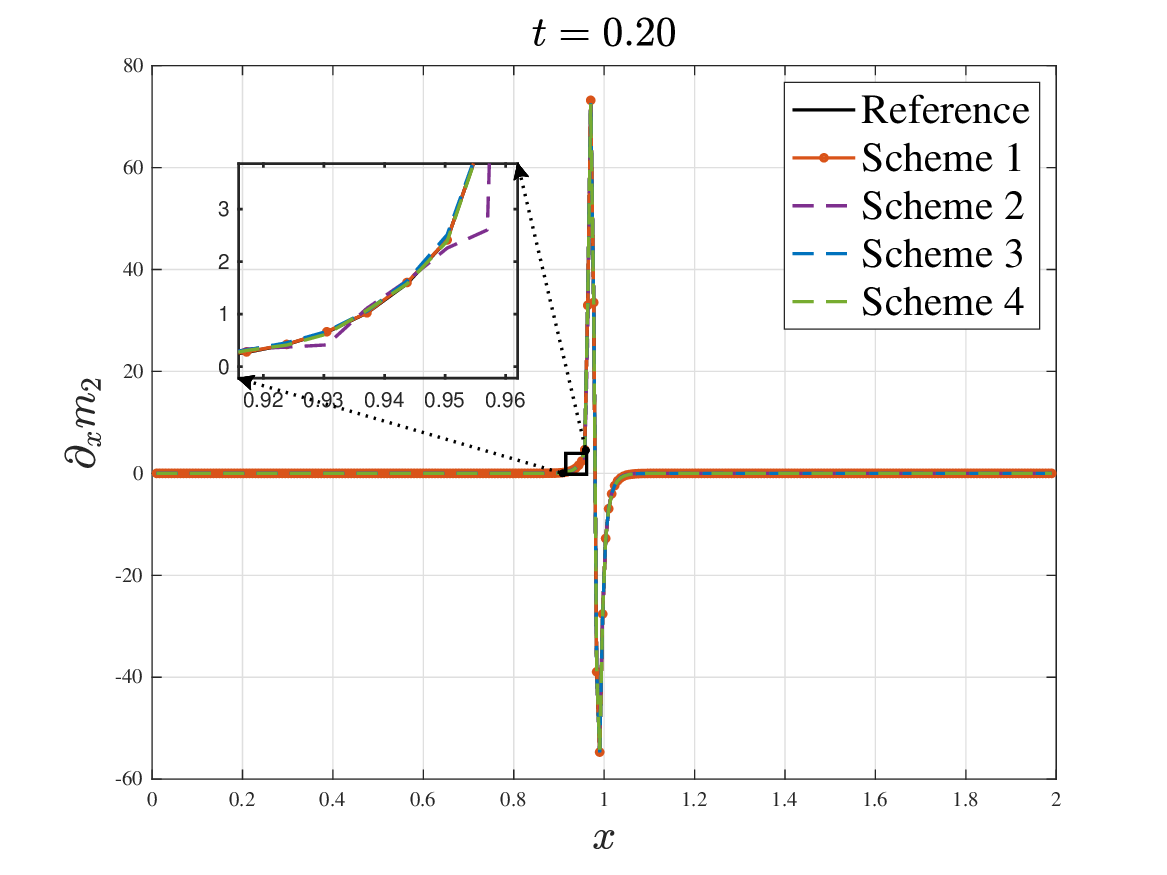} (c)
    \end{minipage}
    \begin{minipage}[b]{0.32\textwidth}
        \centering
\includegraphics[width=\textwidth]{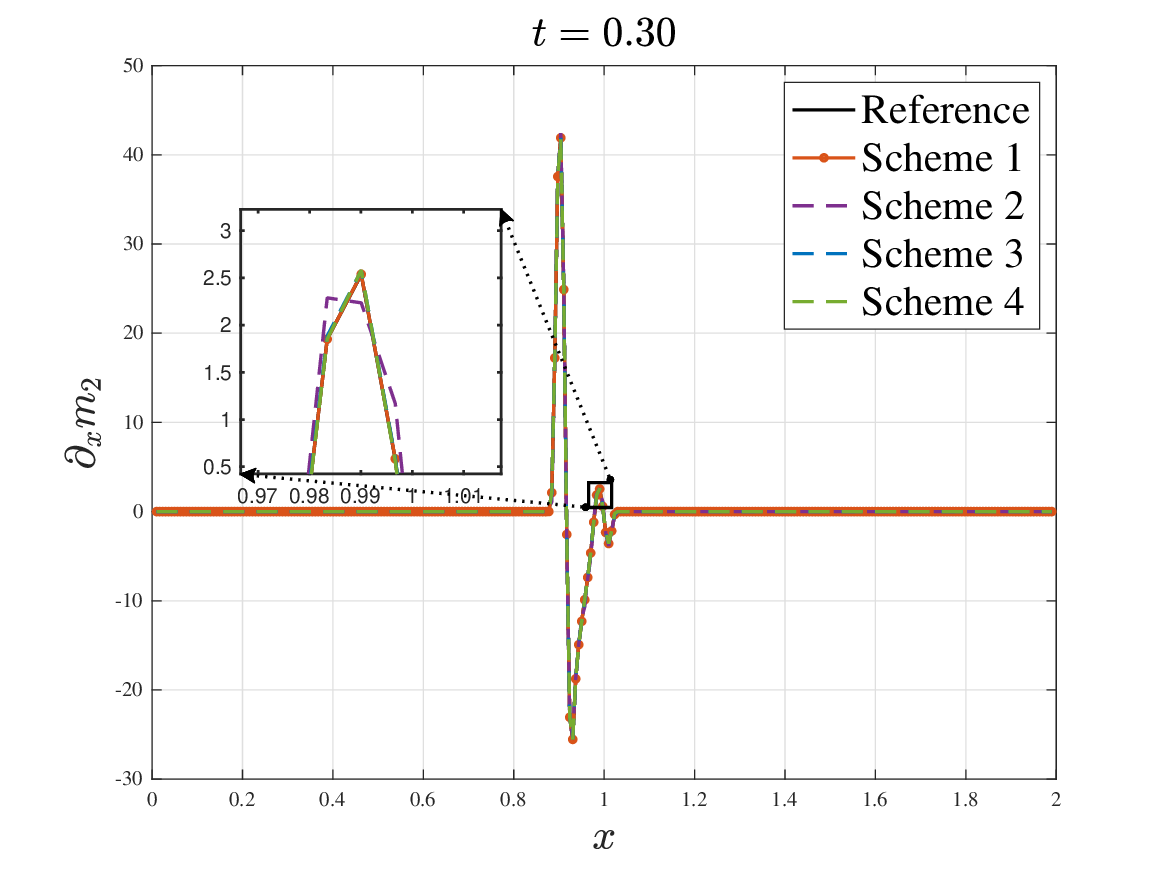} (d)
\end{minipage}
\begin{minipage}[b]{0.32\textwidth}
        \centering
\includegraphics[width=\textwidth]{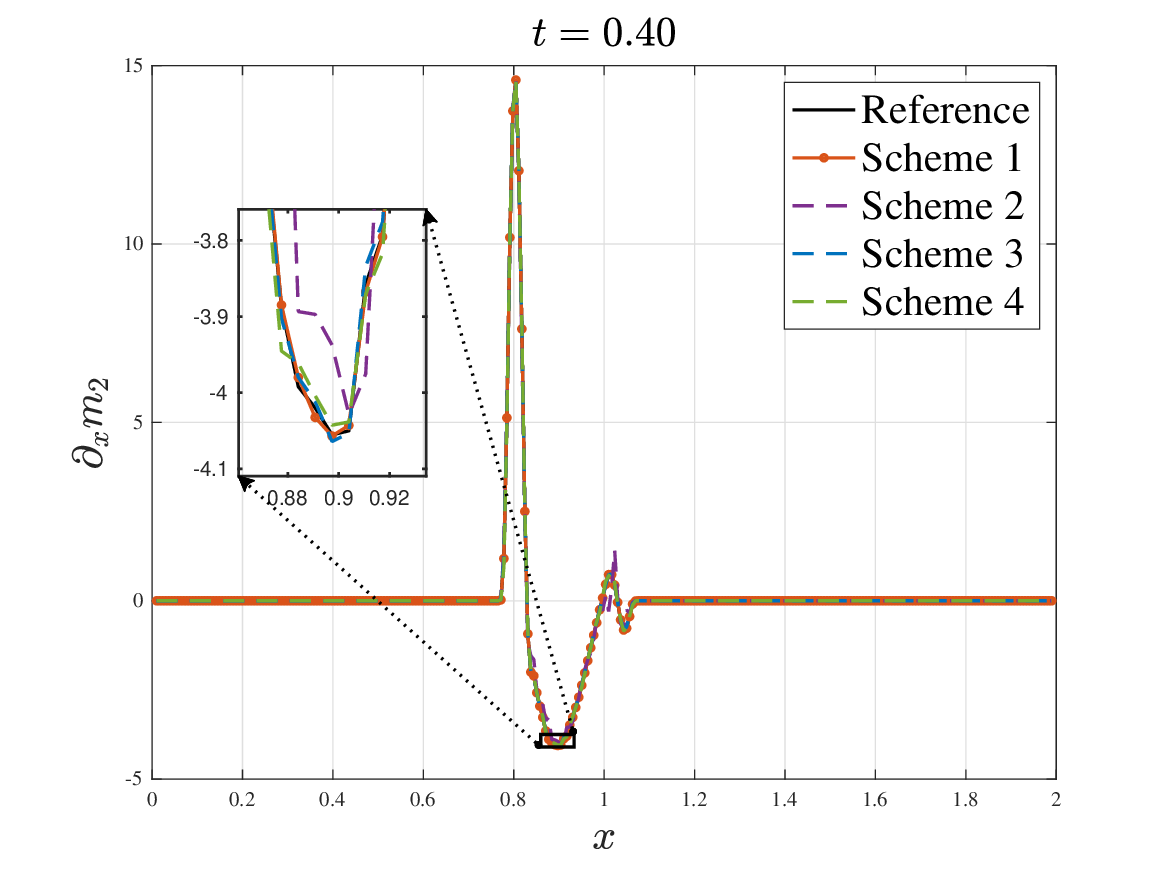} (e)
    \end{minipage}
    \begin{minipage}[b]{0.32\textwidth}
        \centering
        \includegraphics[width=\textwidth]{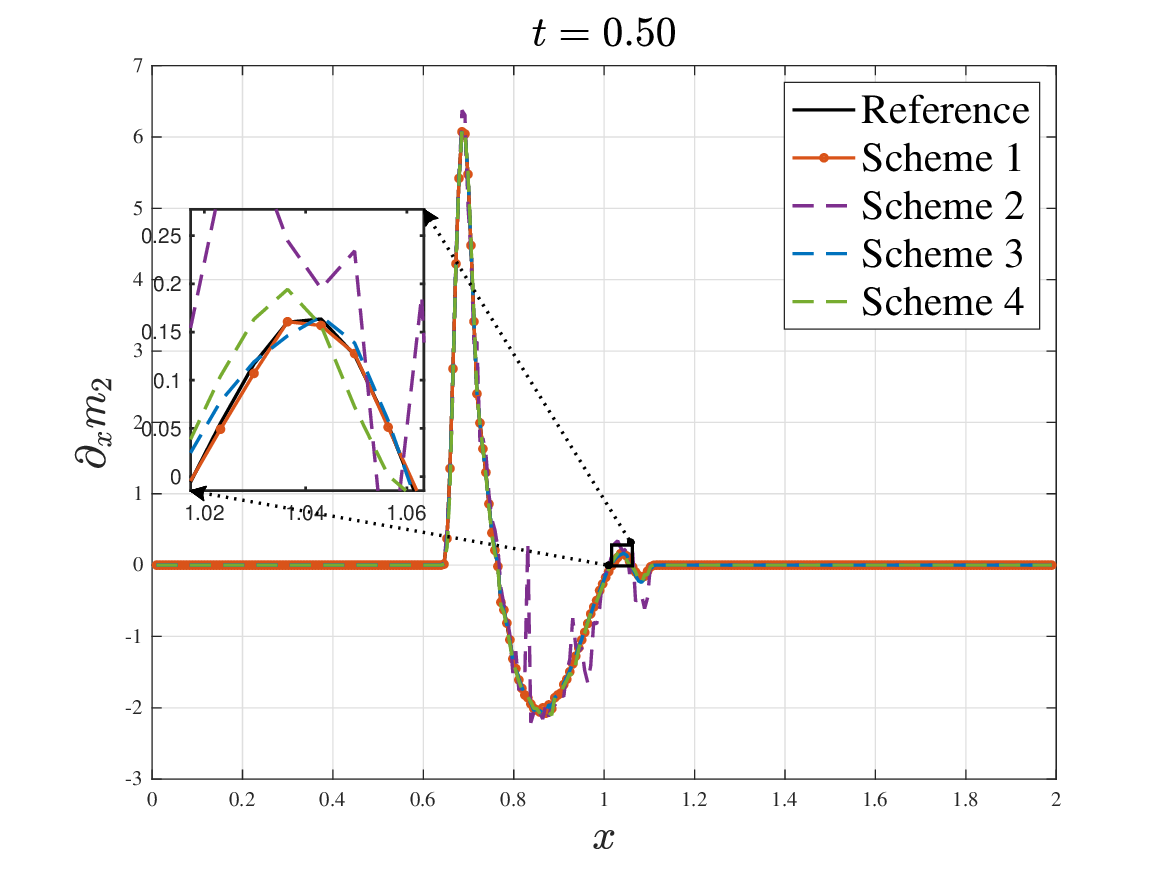} (f)
    \end{minipage}
\caption{Problems I Stage 1 with different $t$. Moment $\partial_x m_2$  for different schemes and reference solutions at different time steps. }
\label{fig:test1_stage1_v1}
\end{figure}

\noindent \textbf{Stage 1.} We first study Stage 1 to compare different schemes with different time steps $T$.
In Figure \ref{fig:test1_stage1_v1}, we plot the reference results for $\partial_x m_2$ obtained by particle in cell method and the approximated solution $\partial_x m_2^{\text{NN},\textit{stage1}}$ generated by different schemes. We observe that the preferred scheme 1 can better approximate ground truth $\partial_x m_2$ under different $T$ and in all regions, while other schemes cannot capture the solutions with high accuracy, especially in some complex regions. 

To qualitatively demonstrate the performance of scheme 1, we report the relative $\ell^2$ errors generated by different schemes. It can be seen that the relative $\ell^2$ error of scheme 1 is much smaller than that of other schemes from Table \ref{tab:test1_stage1}.\\

\begin{table}[htbp]
    \centering
    \vspace{-15pt}
    \begin{small}
            \renewcommand{\multirowsetup}{\centering}
            \setlength{\tabcolsep}{5.5pt}
            \scalebox{0.8}{
                \begin{tabular}{lcccccc}
                    \toprule
                    \multirow{2}{*}{Method} & \multicolumn{6}{c}{time step ($\mathcal{T}$)} \\
                    \cmidrule(lr){2-7}
                    & 0.01 & 0.1 & 0.2 & 0.3 & 0.4 & 0.5 \\
                    \midrule
                    \textbf{Scheme 1} & $\mathbf{7.068 \times 10^{-3}}$ & $\mathbf{1.658 \times 10^{-3}}$ & $\mathbf{1.919 \times 10^{-4}}$ & $\mathbf{1.129 \times 10^{-4}}$ & $\mathbf{1.222 \times 10^{-3}}$ & $\mathbf{2.551 \times 10^{-3}}$ \\
                    Scheme 2 & $3.402 \times 10^{-1}$ & $1.031 \times 10^{-1}$ & $6.894 \times 10^{-2}$ & $2.539 \times 10^{-2}$ & $6.702 \times 10^{-2}$ & $1.849 \times 10^{-1}$ \\
                    Scheme 3 & $3.047 \times 10^{-2}$ & $6.825 \times 10^{-3}$ & $3.015 \times 10^{-3}$ & $1.889 \times 10^{-3}$ & $3.478 \times 10^{-3}$ & $1.228 \times 10^{-2}$ \\
                    Scheme 4 & $8.874 \times 10^{-2}$ & $1.974 \times 10^{-2}$ & $9.273 \times 10^{-3}$ & $2.204 \times 10^{-3}$ & $1.193 \times 10^{-2}$ & $2.920 \times 10^{-2}$ \\
                    \bottomrule
                \end{tabular}
          }
            \caption{Problems I Stage 1. Relative $\ell^2$ error between $\partial_x m_2$ and $\partial_x m_2^{\text{NN},\textit{stage1}}$ with different schemes at different time steps. }
 \label{tab:test1_stage1}
    \end{small}
\end{table}

\begin{figure}[htbp]
    \centering
    \begin{minipage}[b]{0.32\textwidth}
        \centering
        \includegraphics[width=\textwidth]{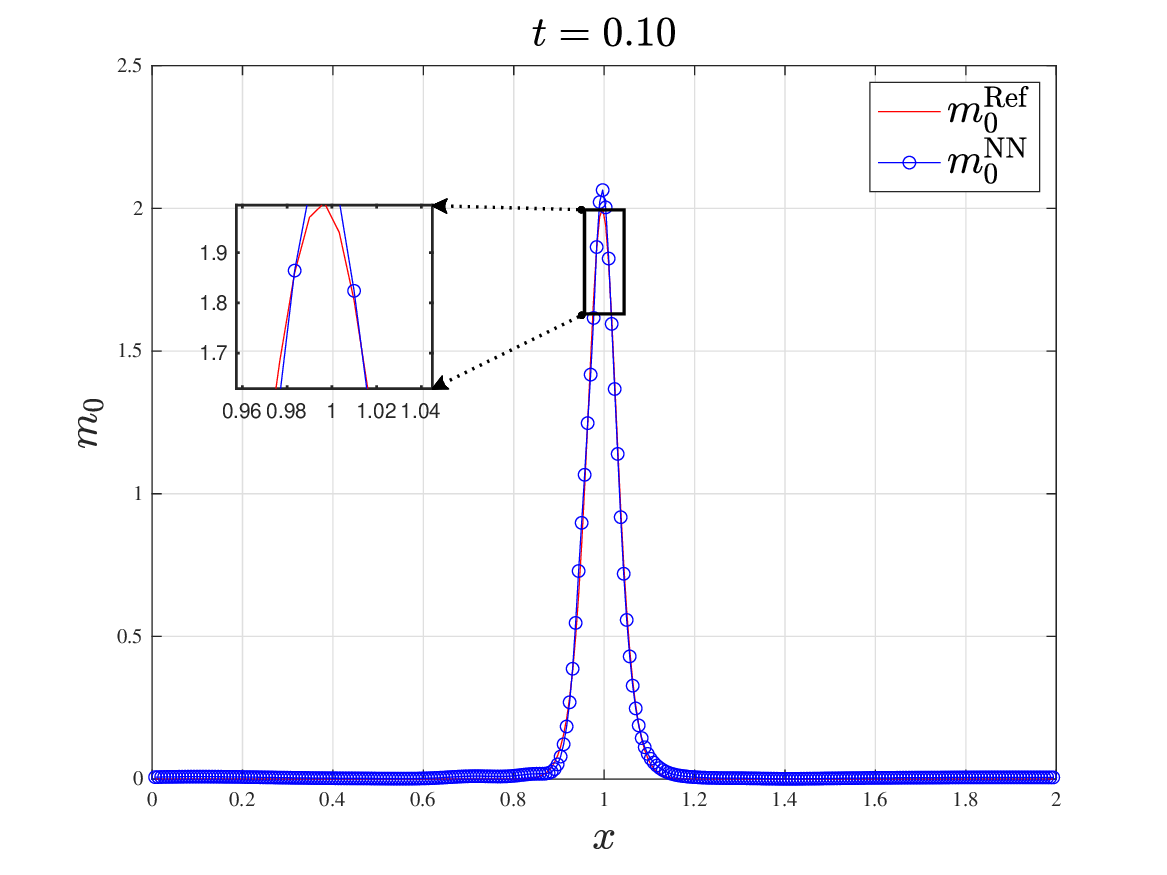} (a)
    \end{minipage}
            \begin{minipage}[b]{0.32\textwidth}
        \centering
        \includegraphics[width=\textwidth]
        {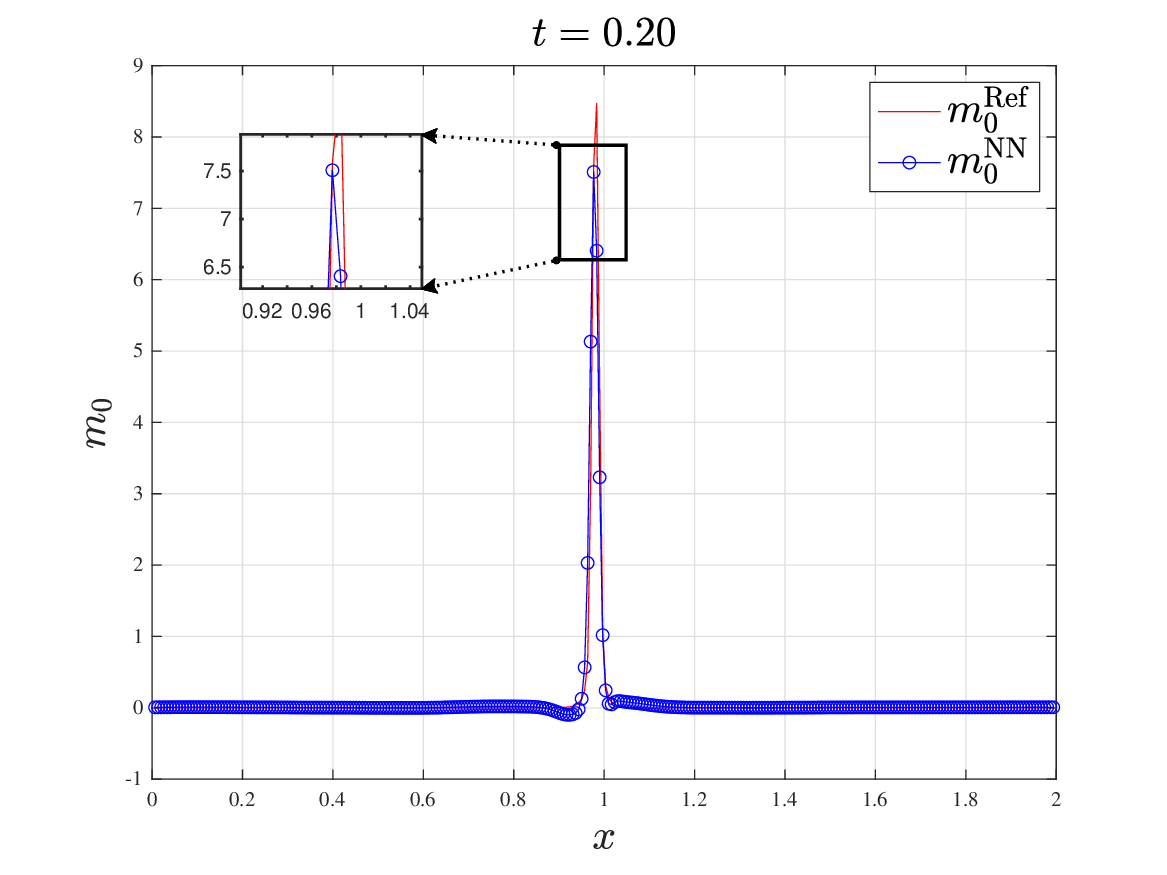} (b)
    \end{minipage}
            \begin{minipage}[b]{0.32\textwidth}
        \centering
        \includegraphics[width=\textwidth]
        {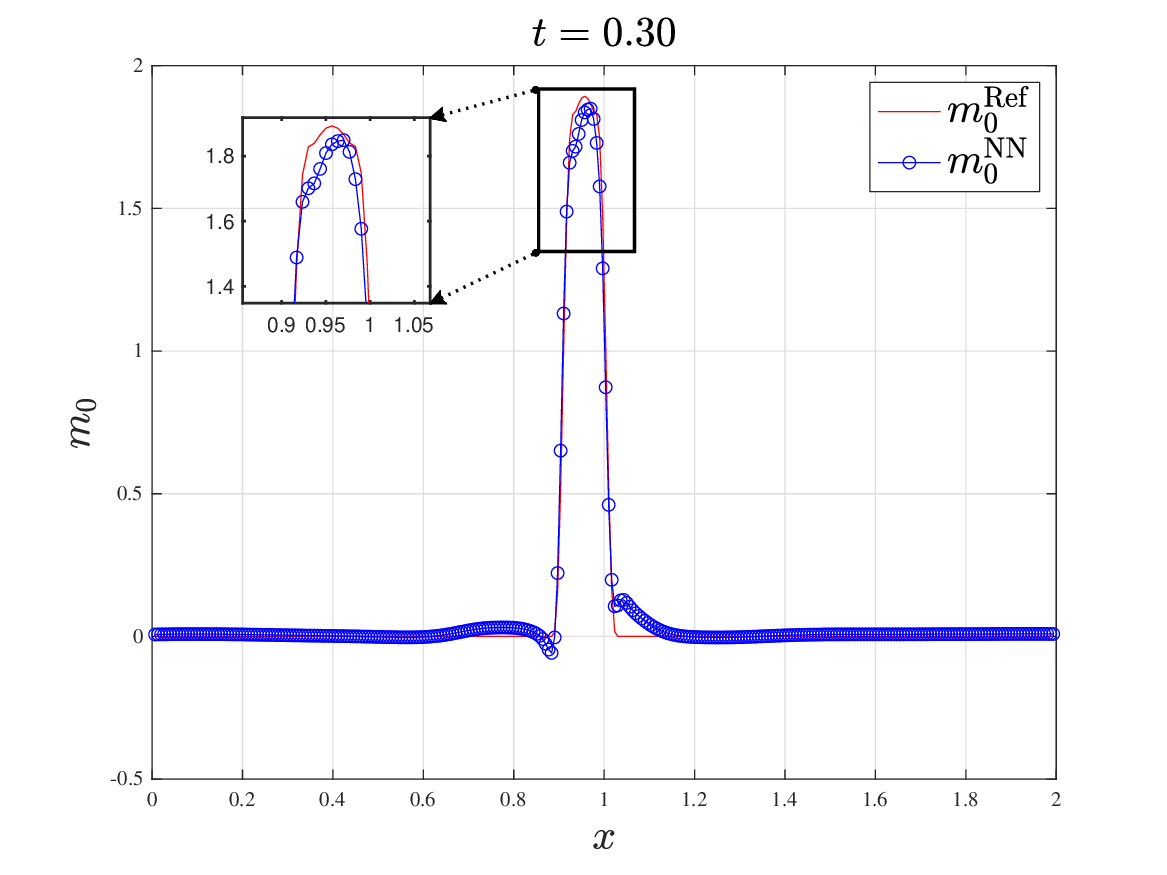}
        (c)
    \end{minipage}
    \begin{minipage}[b]{0.32\textwidth}
        \centering
        \includegraphics[width=\textwidth]
        {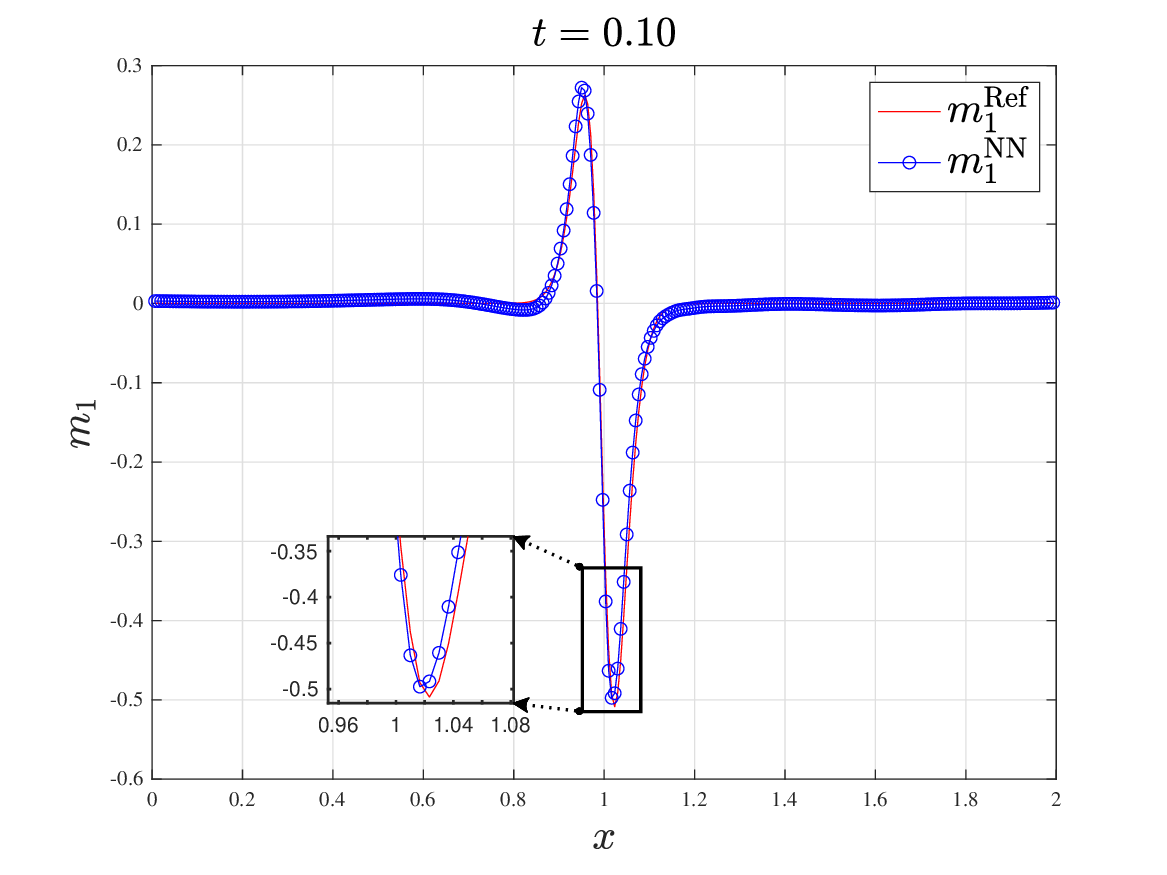} (d)
    \end{minipage}
        \begin{minipage}[b]{0.32\textwidth}
        \centering
        \includegraphics[width=\textwidth]
        {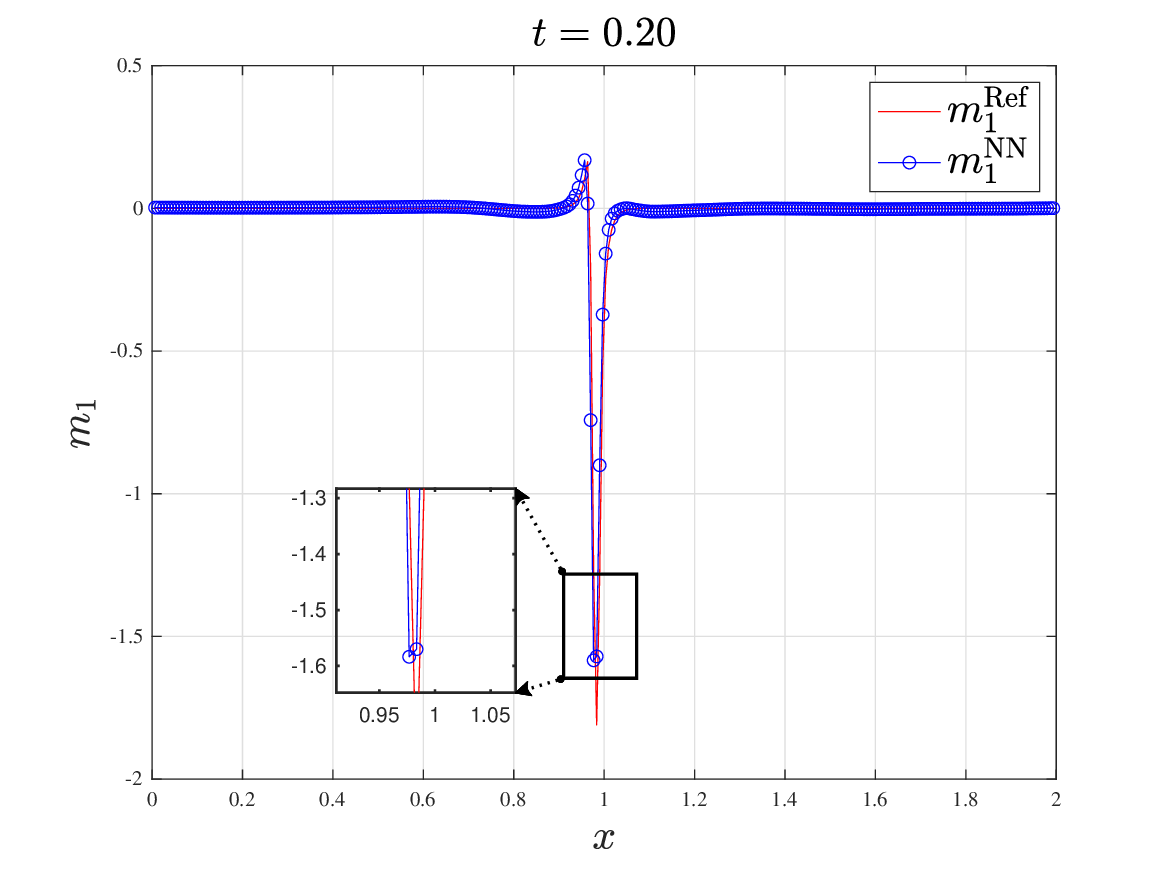} 
        (e)
    \end{minipage}
    \begin{minipage}[b]{0.32\textwidth}
        \centering
        \includegraphics[width=\textwidth]
        {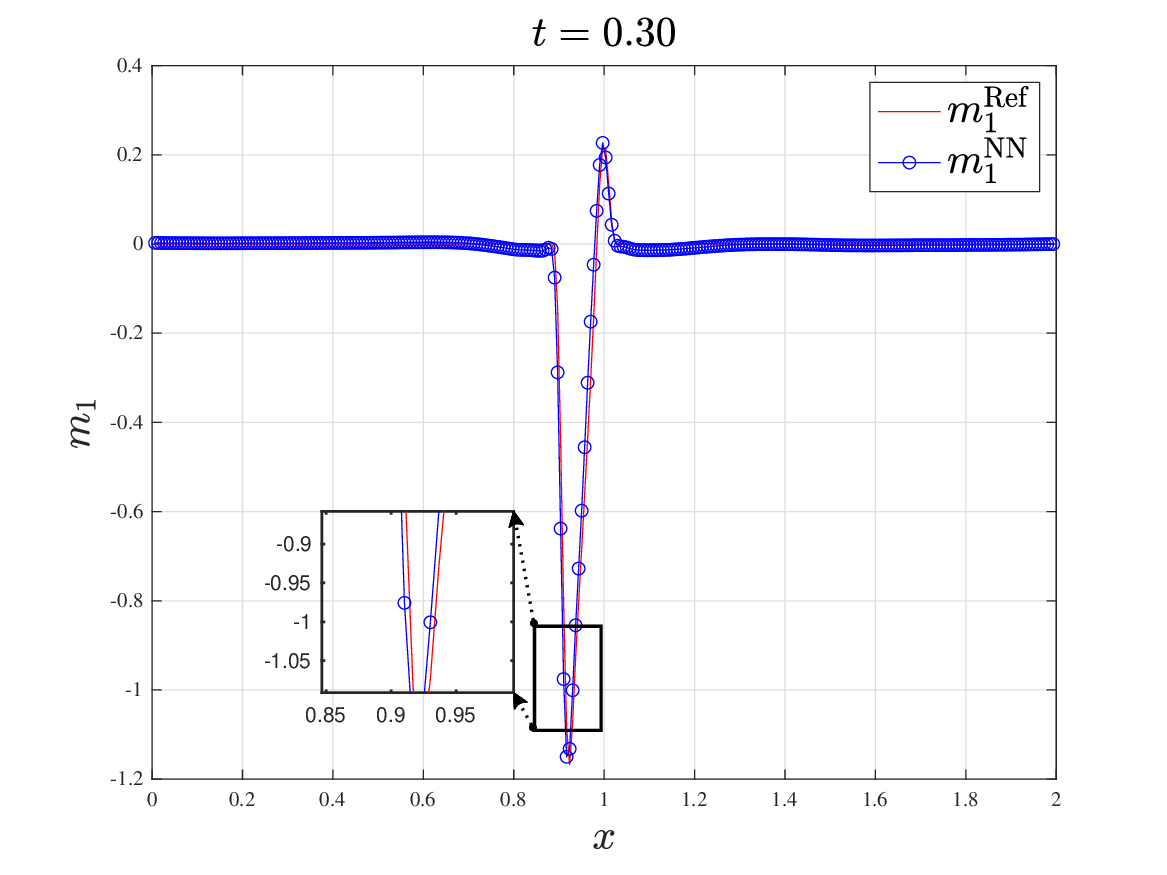}
        (f)
    \end{minipage}
  \caption{Problems I Stage 2 with different $t$. Moments $m_0,m_1$ for our proposed method and reference solutions at different time steps. }
\label{fig:test1_stage2}
\end{figure}

\noindent \textbf{Stage 2.} We further study Stage 2 to compare different moments with different time steps $T$. In Figure \ref{fig:test1_stage2}, we plot the reference results for $m_0, m_1$ and the approximated solutions $m_0^{\text{NN},\textit{stage2}},m_1^{\text{NN},\textit{stage2}}$ obtained by by PINNs at some time steps. We observe that the approximated moments agree with the reference moments under different time steps. The detailed results can be found in Table \ref{tab:test1_m0_m1}.

\begin{table}[htbp]
    \centering
    \begin{small}
            \setlength{\tabcolsep}{5.5pt}
                \begin{tabular}{lccc}
                    \toprule
                    $\mathcal{T}$ & 0.1 & 0.2 & 0.3 \\
                    \midrule
                    $m_0$ & $4.129 \times 10^{-2}$ & $2.757 \times 10^{-1}$ & $8.850 \times 10^{-2}$ \\
                    $m_1$ & $1.009 \times 10^{-1}$ & $3.083 \times 10^{-1}$ & $1.655 \times 10^{-1}$ \\
                    \bottomrule
                \end{tabular}
            \caption{Problems I Stage 2. Relative $\ell^2$ error between the reference moments $m_0$, $m_1$, and the approximated moments $m_0^\text{NN}$, $m_1^\text{NN}$ at different time steps.}
            \label{tab:test1_m0_m1}
    \end{small}
\end{table}

\subsection{Test II: Moment Closure for 1D Problem with Multi-Valued Solutions}
In this test, we consider the problem with multi-valued solutions. Assume the initial data of the Vlasov equation \eqref{Vlasov-eqn} given by 
\begin{equation*} w(t=0,x,v) = \rho_0(x) 
\delta(v - u_0(x)), 
\end{equation*}
where 
\begin{equation*} \rho_0(x) = 1, \quad 
u_0(x) = \chi_{\{x<0\}}-\chi_{\{x>0\}}. 
\end{equation*}
The potential function is given by $\Phi=x^2/2$ (harmonic oscillator). Let the spatial domain be $[x_l, x_r]$ with $x_l=-0.5$ and $x_r=0.5$. The Neumann boundary condition for $m_0$ and $m_1$ is given by $\partial_x m_0(x_{l},t)=\partial_x m_0(x_{r},t) = 0$, $\partial_x m_1(x_{l},t)=\partial_x m_1(x_{r},t) = 0$.
We set $N_x=300$, $\Delta t=5 \times 10^{-3}$ and the final computational time $T=0.2$. The multi-valued solution has two branches, $u_{1,2}(t,x)=-x\tan\left(t\right)\pm\sec\left(t\right)$, see \cite[Example 2]{Jin-Xiantao}.

The main training process for Stage 1 and Stage 2 is similar to that in Test I. For network architecture in Test II, we use a deeper fully connected layer with $10$ layers to better capture shocks in moments.

\vspace{2mm}

\begin{figure}[htbp]
    \centering
    \begin{minipage}[b]{0.46\textwidth}
        \centering
        \includegraphics[width=\textwidth]{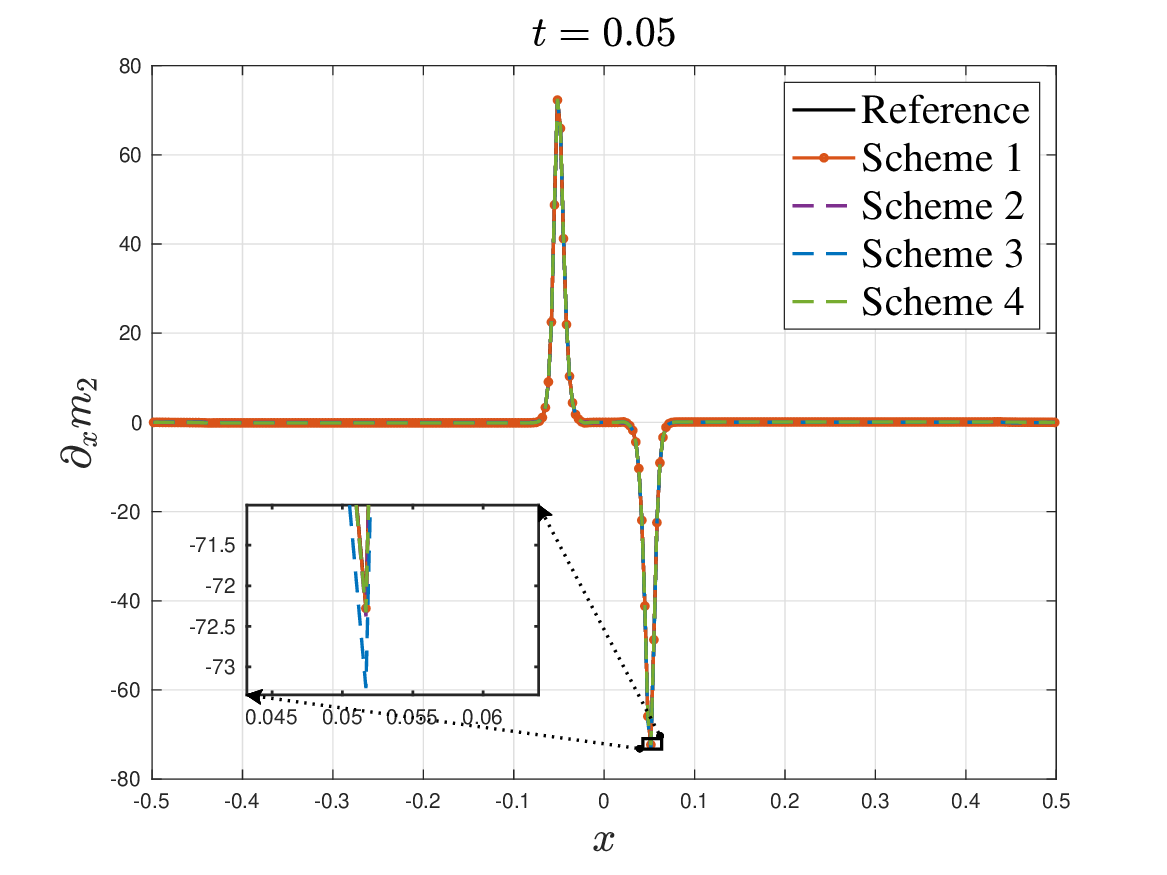} (a)
    \end{minipage}
    \begin{minipage}[b]{0.46\textwidth}
        \centering
        \includegraphics[width=\textwidth]{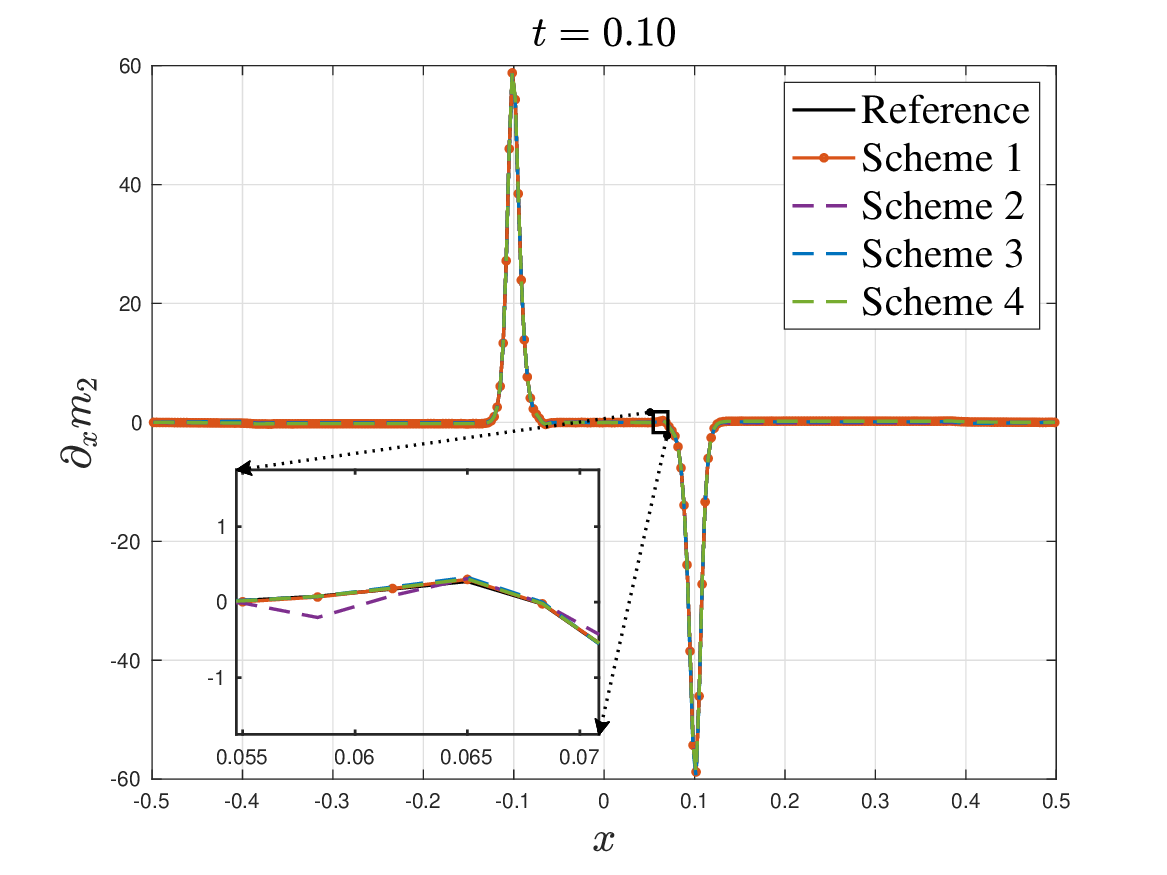} (b)
    \end{minipage}
    \begin{minipage}[b]{0.46\textwidth}
        \centering
        \includegraphics[width=\textwidth]{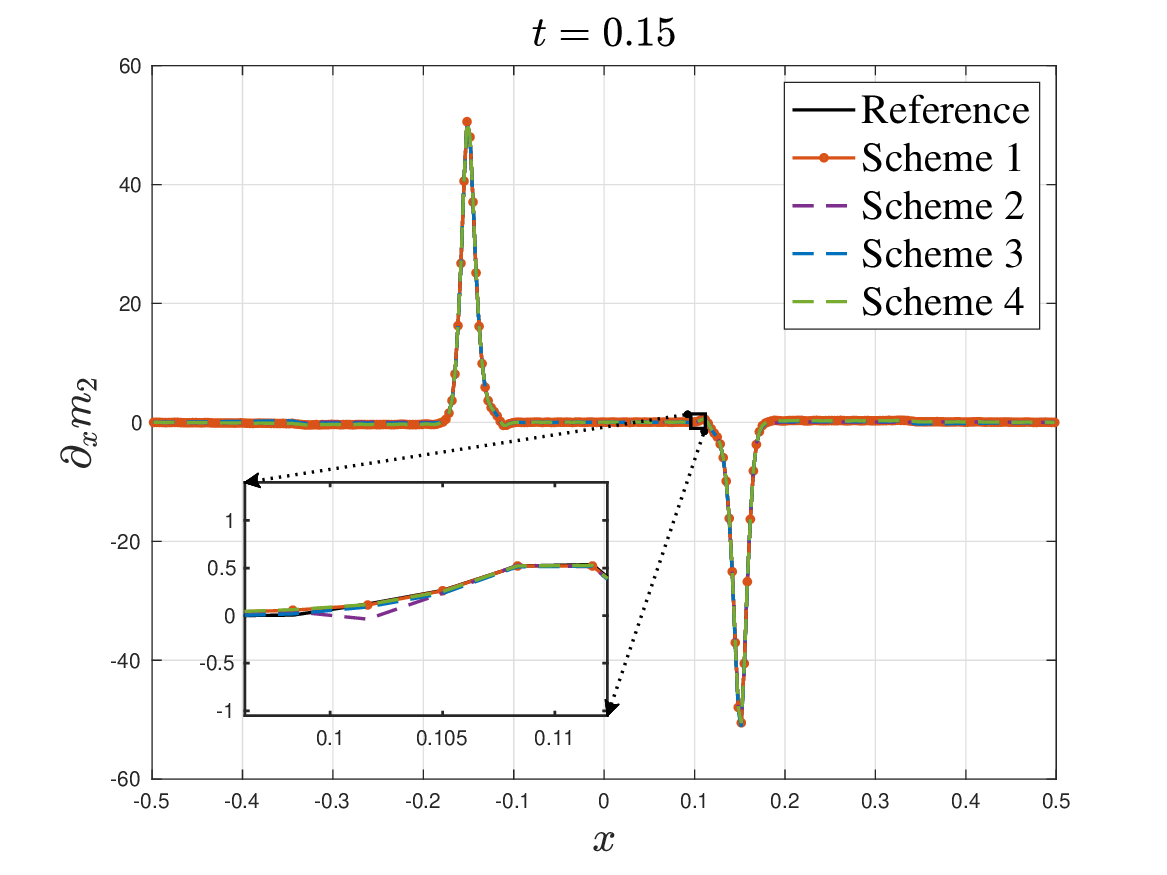} (c)
    \end{minipage}
    \begin{minipage}[b]{0.46\textwidth}
        \centering
        \includegraphics[width=\textwidth]{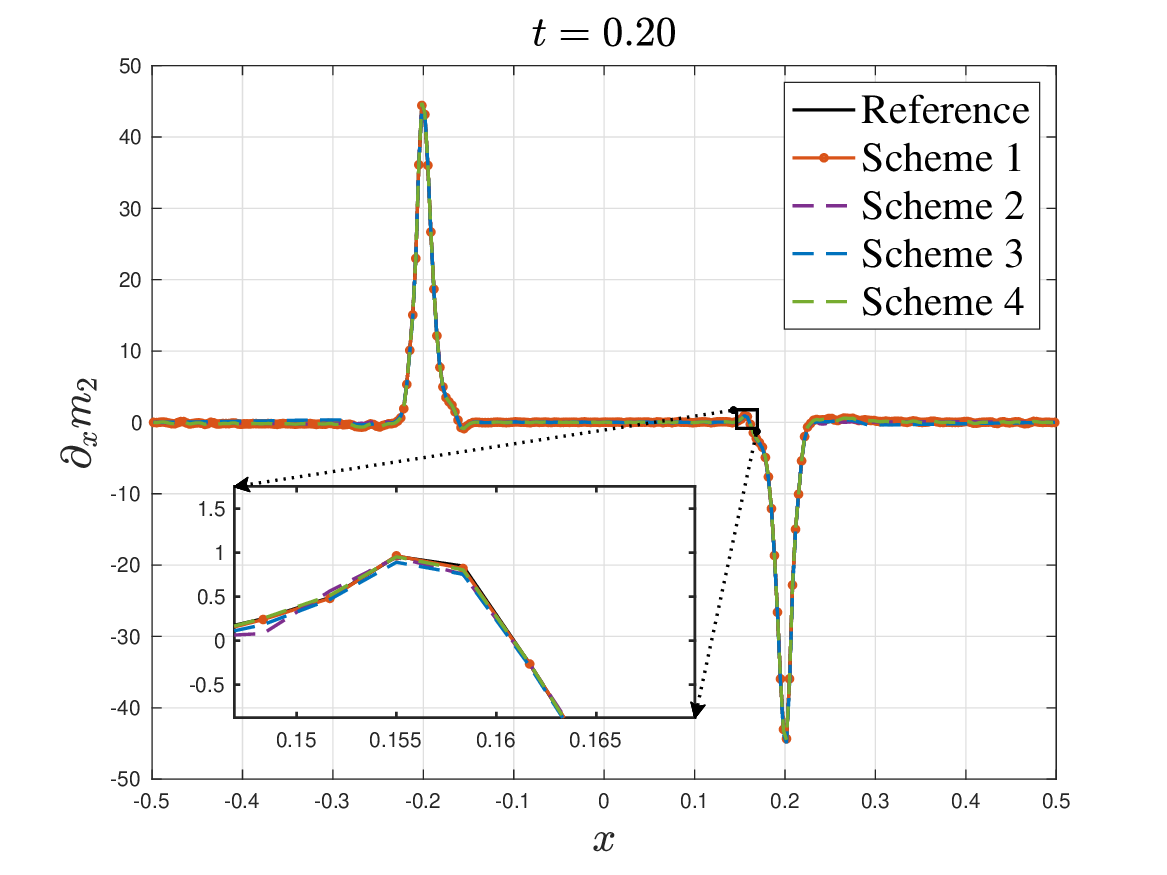} (d)
    \end{minipage}
\caption{Problems II Stage 1 with different $t$. Moment $\partial_x m_2$  for different schemes and reference solutions at time steps $T=0.05, 0.10, 0.15, \text{and}, 0.20$.}
    \label{fig:test2_stage1}
\end{figure}

\noindent \textbf{Stage 1.} We first study Stage 1 to compare different schemes with different time steps T. In Figure \ref{fig:test2_stage1}, we plot the reference results for $\partial_x m_2$ obtained by particle in cell method and the approximated solution $\partial_x m_2^{\text{NN},\textit{stage1}}$ generated by different schemes. It can be seen that the preferred scheme 1 can better approximate the reference moment $\partial_x m_2$ under different $T$ and in all regions, while other schemes cannot capture the solutions with high accuracy, especially in some complex regions.
To qualitatively demonstrate the performance of scheme 1, we report the relative $\ell^2$ errors generated by different schemes. It can be seen that the relative $\ell^2$ error of scheme 1 is much smaller than that of other schemes from  Table \ref{tab:test2_stage1}. \\

\begin{table}[htbp]
    \centering
    \begin{small}
            \renewcommand{\multirowsetup}{\centering}
            \setlength{\tabcolsep}{5.5pt}
            
                \begin{tabular}{lcccc}
                    \toprule
                    \multirow{2}{*}{Method} & \multicolumn{4}{c}{time step ($\mathcal{T}$)} \\
                    \cmidrule(lr){2-5}
                    & 0.05 & 0.10 & 0.15 & 0.20 \\
                    \midrule
                    \textbf{Scheme 1} & $\mathbf{6.288 \times 10^{-4}}$ & $\mathbf{1.254 \times 10^{-3}}$ & $\mathbf{2.371 \times 10^{-3}}$ & $\mathbf{3.401 \times 10^{-3}}$ \\
                    Scheme 2 & $1.043 \times 10^{-2}$ & $1.882 \times 10^{-2}$ & $2.197 \times 10^{-2}$ & $2.397 \times 10^{-2}$ \\
                    Scheme 3 & $1.166 \times 10^{-2}$ & $1.687 \times 10^{-2}$ & $2.472 \times 10^{-2}$ & $3.670 \times 10^{-2}$ \\
                    Scheme 4 & $2.869 \times 10^{-3}$ & $3.511 \times 10^{-3}$ & $3.959 \times 10^{-3}$ & $5.639 \times 10^{-3}$ \\
                    \bottomrule
                \end{tabular}
             \caption{Problems II Stage 1. Relative $\ell^2$ error between $\partial_x m_2$ and $\partial_x m_2^\text{NN}$ with different schemes at different time steps. }
            \label{tab:test2_stage1}
    \end{small}
\end{table} 

\begin{figure}[htbp]
    \centering
    \begin{minipage}[b]{0.32\textwidth}
        \centering
        \includegraphics[width=\textwidth]{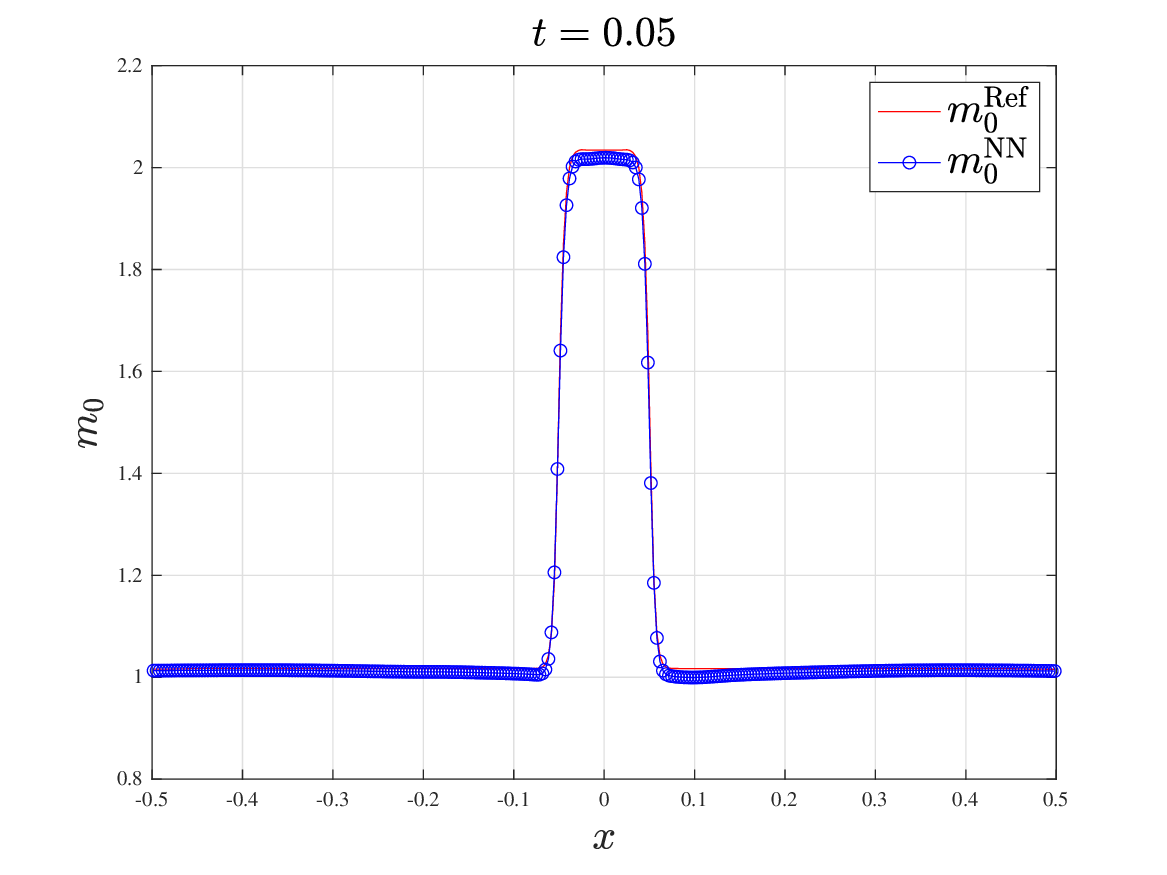} (a)
    \end{minipage}
    \begin{minipage}[b]{0.32\textwidth}
        \centering
        \includegraphics[width=\textwidth]{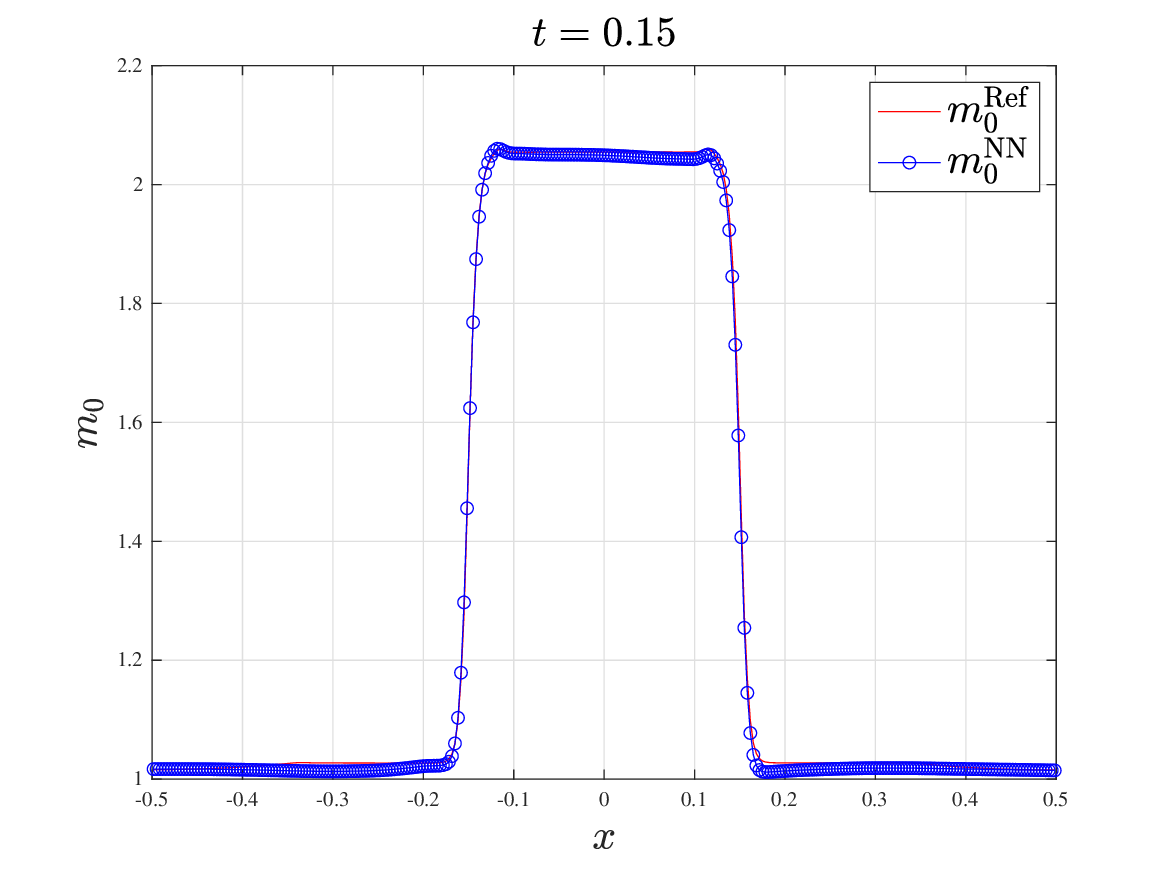} (b)
    \end{minipage}
    \begin{minipage}[b]{0.32\textwidth}
        \centering
        \includegraphics[width=\textwidth]{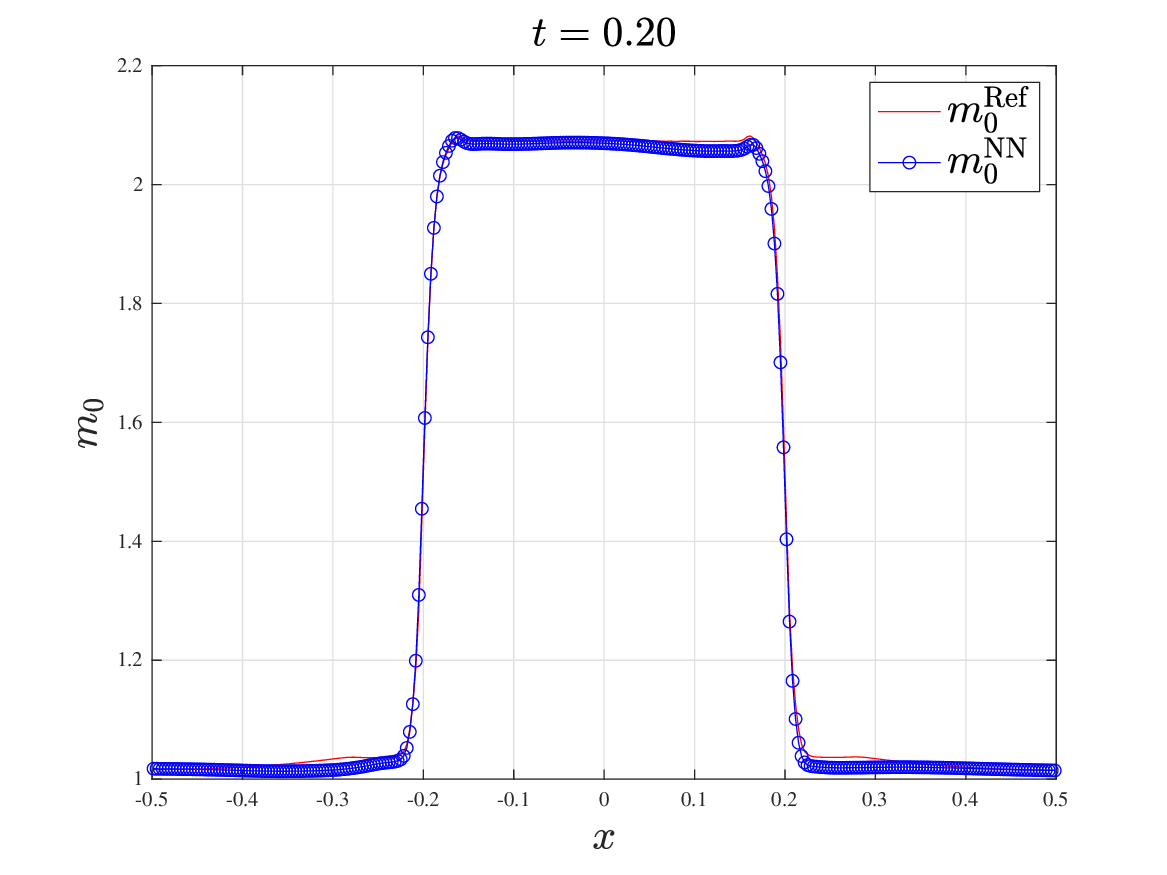} (c)
    \end{minipage}
\caption{Problems II with different $t$. 
$m_0$ for our proposed method and reference solutions at different time steps. }
    \label{fig:test2_stage2_m0}
\end{figure}

\begin{figure}[htbp]
    \centering
    \begin{minipage}[b]{0.32\textwidth}
        \centering
        \includegraphics[width=\textwidth]{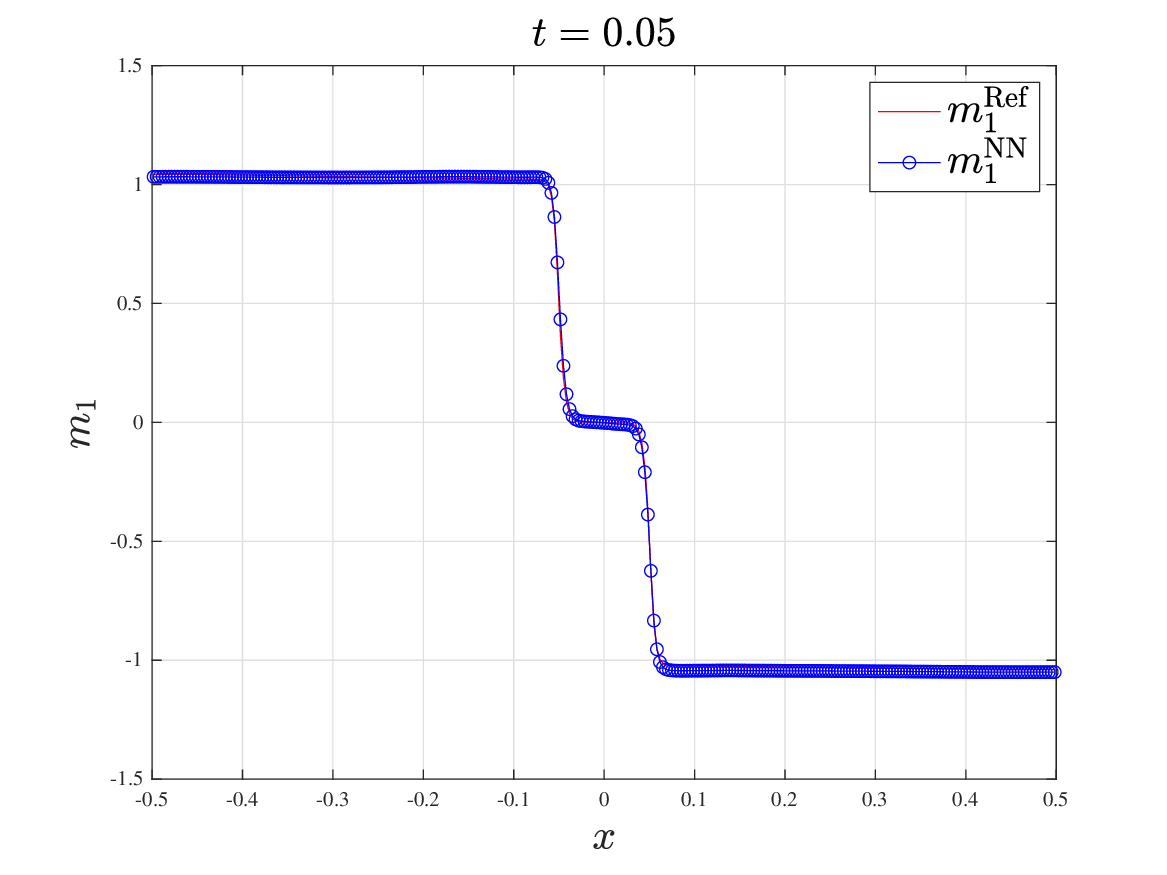} (a)
    \end{minipage}
     \begin{minipage}[b]{0.32\textwidth}
         \centering
         \includegraphics[width=\textwidth]{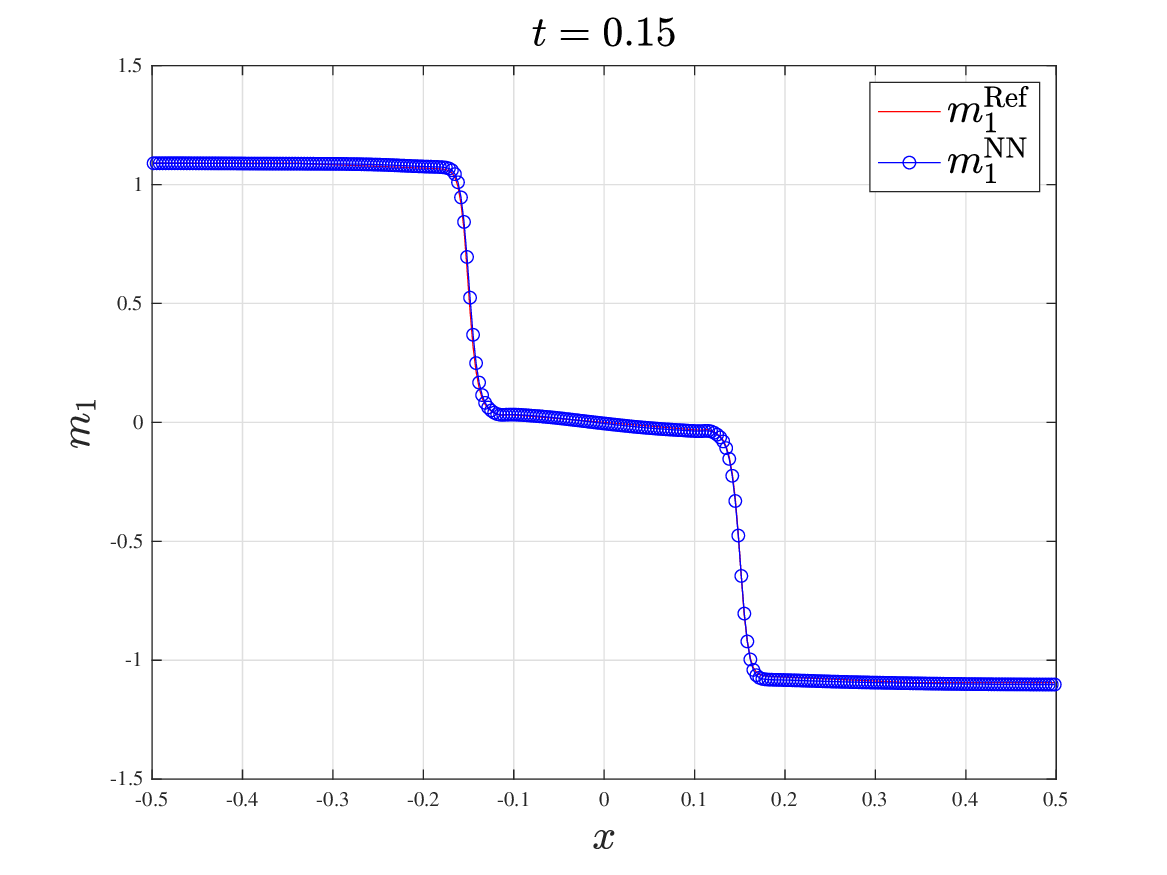} (b)
     \end{minipage}
    \begin{minipage}[b]{0.32\textwidth}
        \centering
        \includegraphics[width=\textwidth]{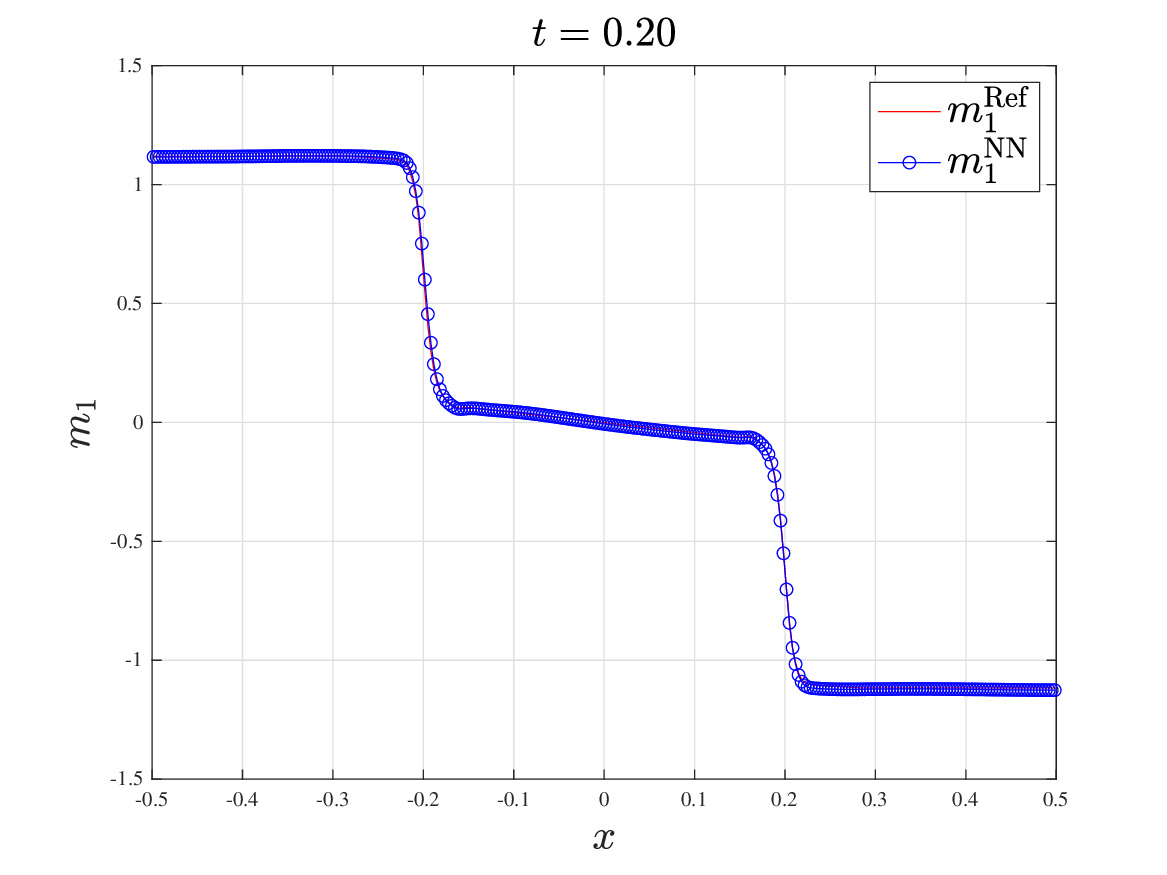} (c)
    \end{minipage}
\caption{Problems II with different $t$. 
$m_1$ for our proposed method and reference solutions at different time steps.  }
    \label{fig:test2_stage2_m1}
\end{figure}

\noindent \textbf{Stage 2.} We further study Stage 2 to compare different moments with different time steps $T$. In Figure \ref{fig:test2_stage2_m0} and Figure \ref{fig:test2_stage2_m1}, we plot the reference results for $m_0, m_1$ and approximated  $m_0^{\text{NN},\textit{stage2}}$,\, $m_1^{\text{NN},\textit{stage2}}$ obtained by by PINNs. We can find from the figures that the moments $m_0^{\text{NN}}$ and $m_1^{\text{NN}}$ match real moments $m_0$ and $m_1$ under different time steps very well. The detailed results can be found in Table \ref{tab:test2_m0_m1}. 

\begin{table}[htbp]
    \centering
    \begin{small}
            \setlength{\tabcolsep}{5.5pt}
   \begin{tabular}{lcccc}
                    \toprule
                    $\mathcal{T}$ & 0.05 & 0.10 & 0.15 & 0.20 \\
                    \midrule
                    $m_0$ & $8.892 \times 10^{-2}$ & $7.109 \times 10^{-3}$ & $7.123 \times 10^{-3}$ & $7.344 \times 10^{-3}$ \\
                    $m_1$ & $1.311 \times 10^{-2}$ & $1.151 \times 10^{-2}$ & $9.865 \times 10^{-3}$ & $9.411 \times 10^{-3}$ \\
                    \bottomrule
                \end{tabular}
             \caption{Problems II Stage 2. Relative $\ell^2$ error between the reference moments $m_0$, $m_1$, and the approximated moments $m_0^\text{NN}$, $m_1^\text{NN}$ at different time steps.}
            \label{tab:test2_m0_m1}
    \end{small}
\end{table}

\subsection{Test III: Moment Closure for 2D Problem} 
\label{Test3}
In this test, we consider the $2D$ problem with multi-valued solutions. Assume the initial data of the Vlasov equation \eqref{Vlasov-eqn} given by 
\begin{equation} 
w(x_1,x_2,v,\xi,t=0) = \rho_0
\delta(v - u_0(x_1)) \delta(\xi - u_0(x_2)), 
\end{equation}
where 
\begin{equation} \rho_0 = 1, \quad 
u_0(x_1) = \chi_{\{x_1<0\}}-\chi_{\{x_1>0\}},  
\end{equation}
and            
$$ u_0(x_2) = \chi_{\{x_2<0\}}-\chi_{\{x_2>0\}}. $$
The potential function is given by $\Phi(x_1, x_2) =\dfrac{x_1^2+x_2^2}{2}$ (harmonic oscillator). Let the spatial domain be $[-0.5,0.5] \times [-0.5,0.5]$. The Neumann boundary condition for $m_0$ and $m_1$ is given by $\partial_{x_1} m_0(x_1^l,x_2,t)=\partial_{x_1} m_0(x_1^{r},x_2,t) = \partial_{x_2} m_0(x_1,x_2^l,t) = \partial_{x_2} m_0(x_1,x_2^r,0) = 0$, $\partial_{x_1} m_1(x_1^l,x_2,t)=\partial_{x_1} m_1(x_1^{r},x_2,t) = \partial_{x_2} m_1(x_1,x_2^l,t) = \partial_{x_2} m_1(x_1,x_2^r,t) = 0$.
We set $N_{x_1}=N_{x_2}=40$, $\Delta t=5 \times 10^{-3}$ and the final computational time $T=0.1$. 

\vspace{2mm}

\noindent \textbf{Formulation.} Consider $m_0$, $m_{11}$, $m_{12}$, and $m_2$ in $2D$ case, we have
\begin{equation}
    \begin{cases}
        m_0 =  \int_{v_1} \int_{v_2} w(t,\bx,v_1,v_2) d v_2 d v_1,\\[2pt]
        m_{11} =  \int_{v_1} \int_{v_2} w(t,\bx,v_1,v_2) v_1 d v_2 d v_1, \\m_{12} =  \int_{v_1} \int_{v_2} w(t,\bx,v_1,v_2) v_2 d v_2 d v_1,\\[2pt]
        m_2 = \int_{v_1} \int_{v_2} w(t,\bx,v_1,v_2) (v_1^2 + v_2^2) d v_2 d v_1,\\[2pt]
        m_{21} = \int_{v_1} \int_{v_2} w(t,\bx,v_1,v_2) v_1^2 d v_2 d v_1,\\[2pt]
        m_{22} = \int_{v_1} \int_{v_2} w(t,\bx,v_1,v_2) v_2^2 d v_2 d v_1.
    \end{cases}
\end{equation}
Consider equation \eqref{Vlasov-eqn} in $2D$ case, 
\begin{equation}\label{Vlasov-eqn2D}
   \partial_t w + \boldsymbol{v} \cdot \nabla_{\bx} w - \nabla_{\bx} \Phi \cdot \nabla_{\boldsymbol{v}} w=0. 
\end{equation}

We derive the moment system for the $2D$ case, with details shown in the Appendix:  
\begin{equation} \label{test3_system}
\begin{cases} 
       \partial_t m_0 + \partial_{x_1} m_{11} + \partial_{x_2} m_{12} = \int_{v_1} \int_{v_2} (x_1 \partial_{v_1} w + x_2 \partial_{v_2} w) d v_2 d v_1,\\[2pt]
        \partial_t m_{11} + \int_{v_1} \int_{v_2} (v_1^2 \partial_{x_1} w  +  v_1 v_2 \partial_{x_2} w) +  v_1 (x_1 \partial_{v_1} w + x_2 \partial_{v_2} w) d v_2 d v_1=0,\\[2pt]
        \partial_t m_{12} + \int_{v_1} \int_{v_2} (v_1 v_2 \partial_{x_1} w  +  v_2^2 \partial_{x_2} w) - v_2 (x_1 \partial_{v_1} w + x_2 \partial_{v_2} w) d v_2 d v_1=0. 
\end{cases} 
\end{equation}
When we employ the PINNs on the moment system, the system \eqref{2-2moment} changes to the following form:
\begin{equation} \notag
\hspace{-1mm}\begin{cases} 
       \partial_t m_0^{\text{NN},\textit{stage2}} + \partial_{x_1} m_{11}^{\text{NN},\textit{stage2}} + \partial_{x_2} m_{12}^{\text{NN},\textit{stage2}} = 0, \\[2pt] 
        \partial_t m_{11}^{\text{NN},\textit{stage2}} + \partial_{x_1} m_{21}^{\text{NN},\textit{stage1}} +  \Phi_{x_1} m_{11}^{\text{NN},\textit{stage2}} + \int_{v_1} \int_{v_2} (  v_1 v_2 \partial_{x_2} w + v_1 x_2 \partial_{v_2} w )  d v_2 d v_1=0,\\[2pt]
        \partial_t m_{12}^{\text{NN},\textit{stage2}} + \partial_{x_2} m_{22}^{\text{NN},\textit{stage1}} - \Phi_{x_2} m_{12}^{\text{NN},\textit{stage2}} +\int_{v_1} \int_{v_2}( v_1 v_2 \partial_{x_1} w   - v_2 x_1 \partial_{v_1} w ) d v_2 d v_1=0,
\end{cases} 
\end{equation}
where $\bx = (x_1,x_2)$ is a $2D$ variable.

In Stage 1, we consider learning $\partial_{x_1} m_{21}^{\text{NN}, \textit{stage}_1}$ and  $\partial_{x_2} m_{22}^{\text{NN}, \textit{stage}_1}$ in the following way, which is different from Test I and Test II. For simplicity, we can leverage two different neural networks $\mathcal{NN}_1$ and $\mathcal{NN}_2$ with similar structures to learn $\partial_{x_1} m_{21}^{\text{NN}, \textit{stage}_1}$ and  $\partial_{x_2} m_{22}^{\text{NN}, \textit{stage}_1}$, respectively. The schemes in Test III are also different from those in Test I and Test II. For scheme 1, we consider
\begin{equation}
    \begin{cases}
        \partial_{x_1} m_{21}^{\text{NN}, \textit{stage}_1} = \mathcal{NN}_1(m_0, \partial_{x_1} m_0, \partial_{x_2} m_0, m_{11}, \partial_{x_1} m_{11}, \partial_{x_2} m_{11}, m_{12}, \partial_{x_1} m_{12}, \partial_{x_2} m_{12}),\\[2pt]
        \partial_{x_2} m_{22}^{\text{NN}, \textit{stage}_1} = \mathcal{NN}_2(m_0, \partial_{x_1} m_0, \partial_{x_2} m_0, m_{11}, \partial_{x_1} m_{11}, \partial_{x_2} m_{11}, m_{12}, \partial_{x_1} m_{12}, \partial_{x_2} m_{12}),
    \end{cases}
\end{equation}
More details of other schemes can be found in Appendix \ref{scheme_appen}.

\vspace{2mm}

\begin{figure}[htbp] 
    \centering
    \begin{minipage}[b]{0.49\textwidth}
        \centering
        \includegraphics[width=\textwidth]{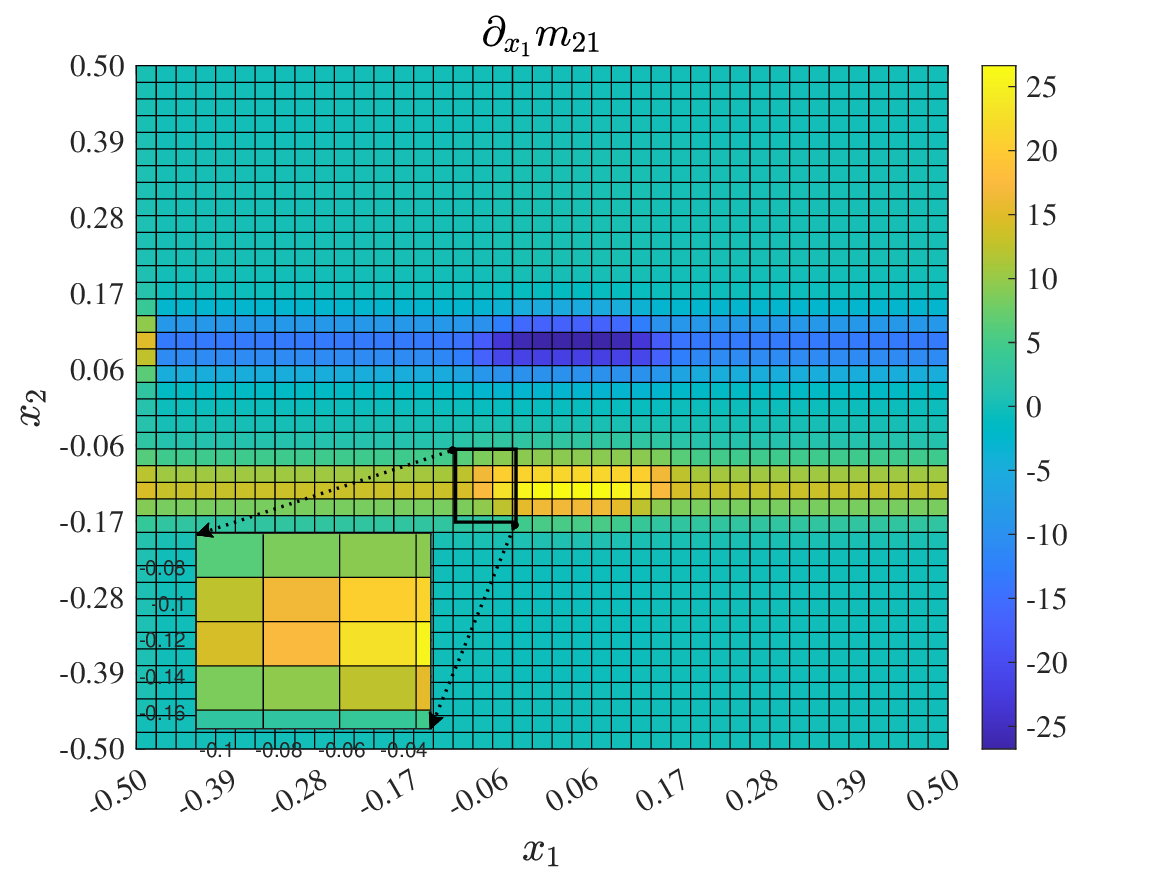}
    \end{minipage}
    \begin{minipage}[b]{0.49\textwidth}
        \centering
        \includegraphics[width=\textwidth]{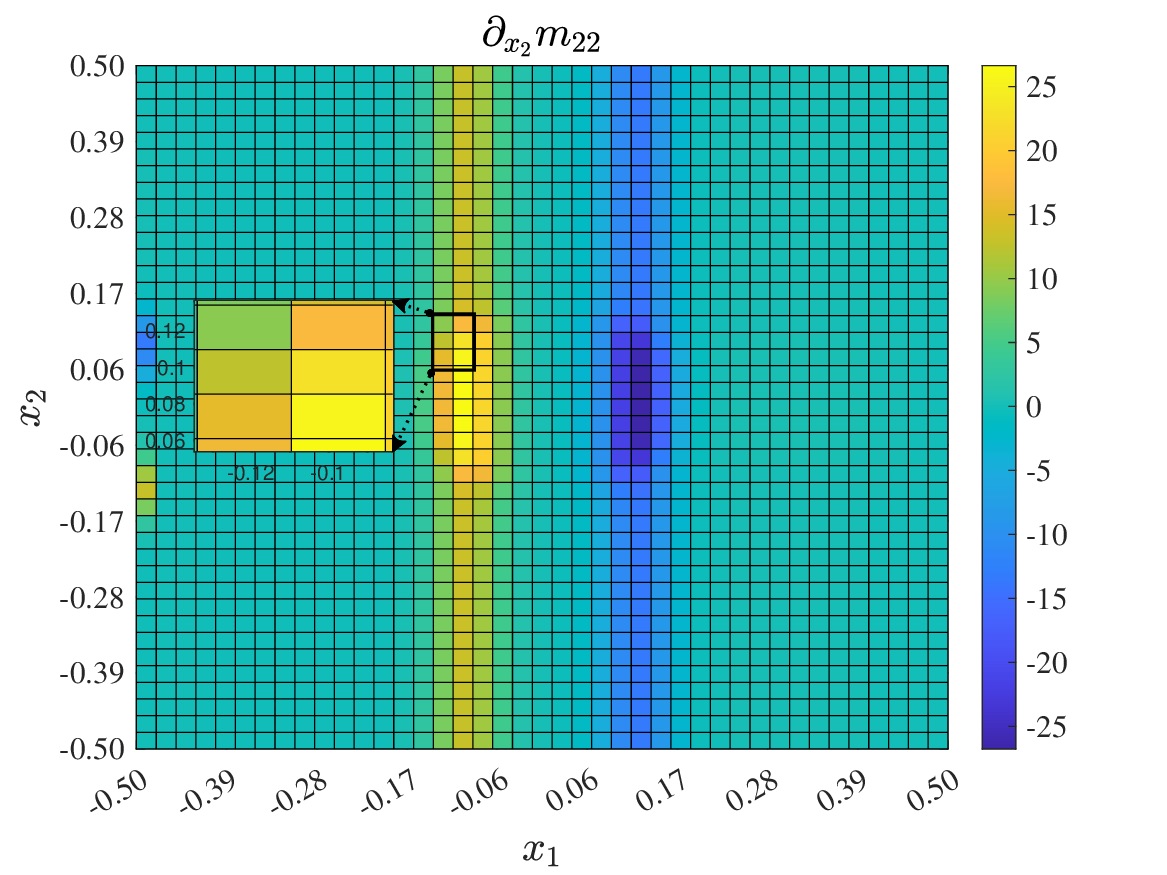}
    \end{minipage}
    \caption{Problems III Stage 1. Reference moments $\partial_{x_1} m_{21}$ and $\partial_{x_2} m_{22}$ solutions at final time $T=0.1$.}
    \label{fig:test3_stage1_ref_v2}  
\end{figure}

\begin{figure}[htbp]
    \centering
    \begin{minipage}[b]{0.49\textwidth}
        \centering
        \includegraphics[width=\textwidth]{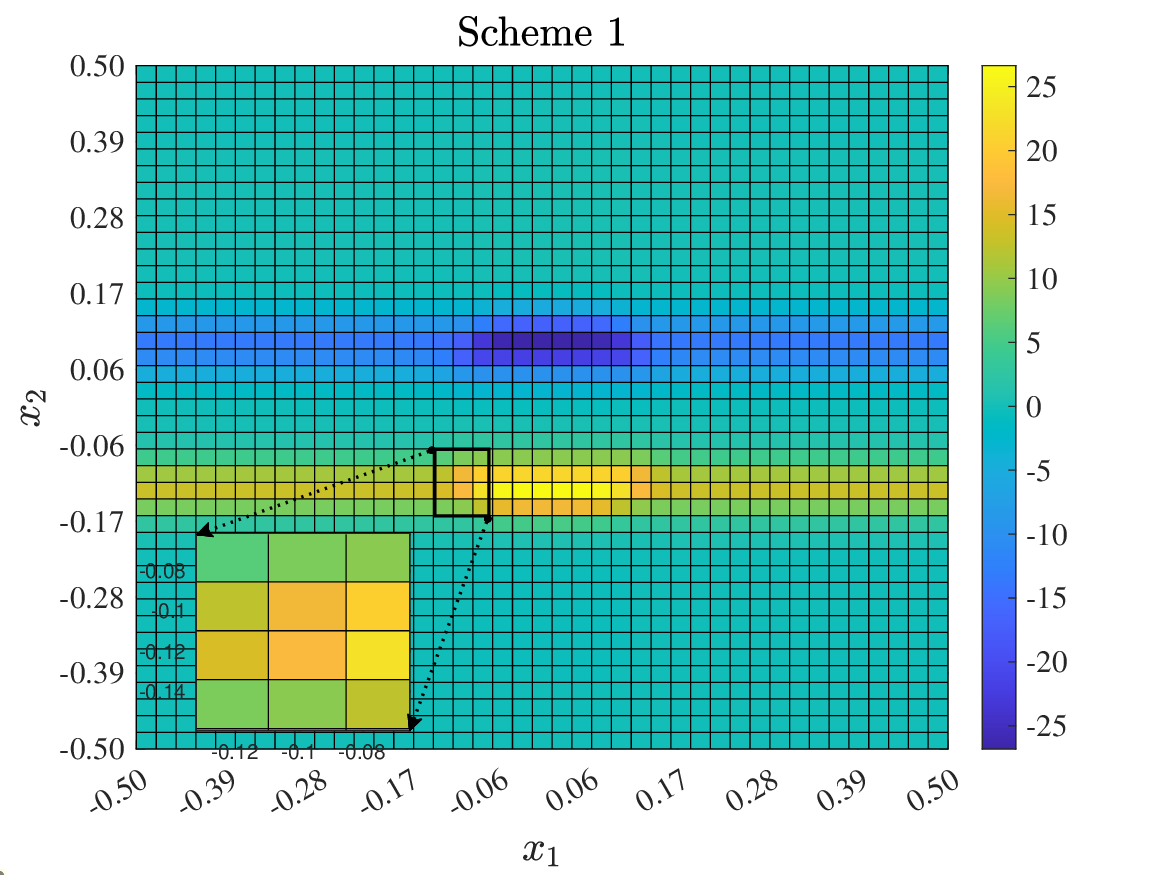} 
    \end{minipage}
    \begin{minipage}[b]{0.49\textwidth}
        \centering
        \includegraphics[width=\textwidth]{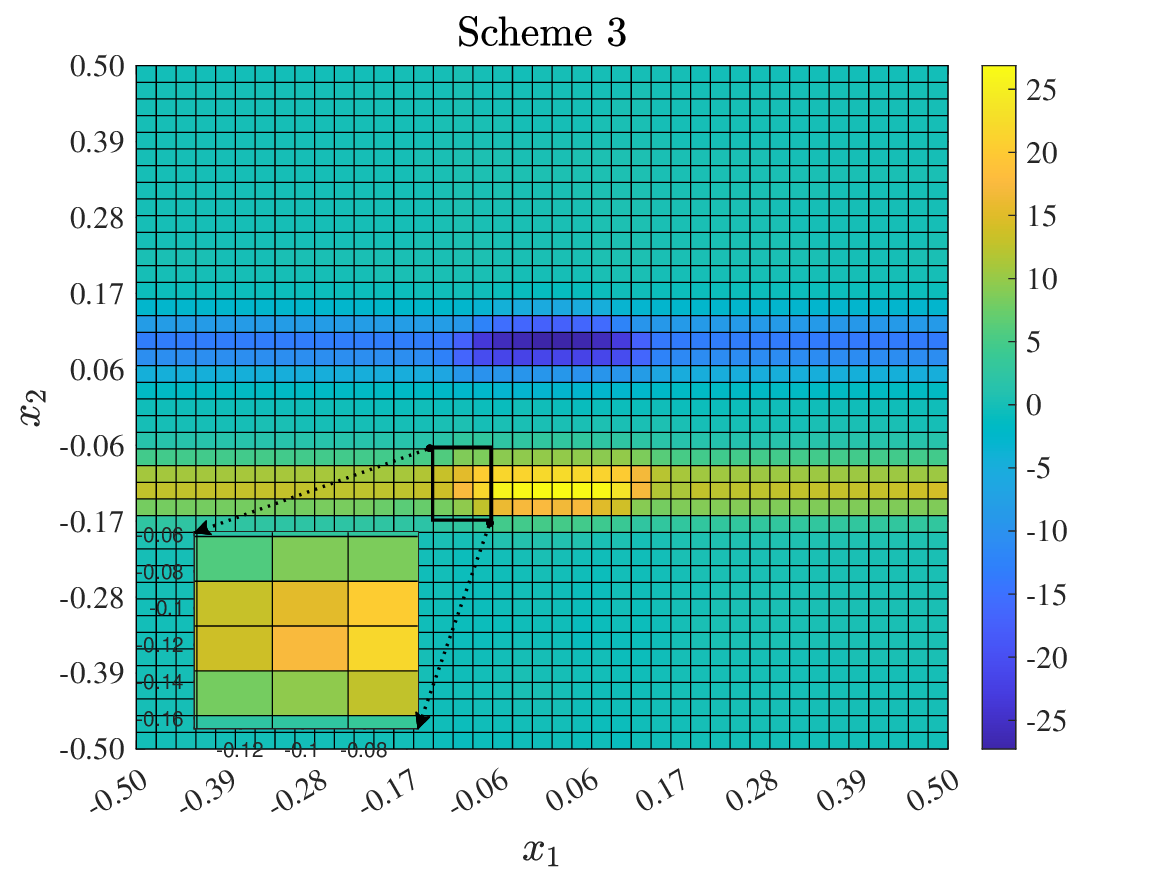} 
    \end{minipage}
\caption{Problems III Stage 1. Moment $\partial_{x_1} m_{21}^{\text{NN}}$ with different schemes at final time $T=0.1$.}
    \label{fig:test3_stage1_m21_v2}
\end{figure}

\begin{figure}[htbp]
    \centering
    \begin{minipage}[b]{0.49\textwidth}
        \centering
        \includegraphics[width=\textwidth]{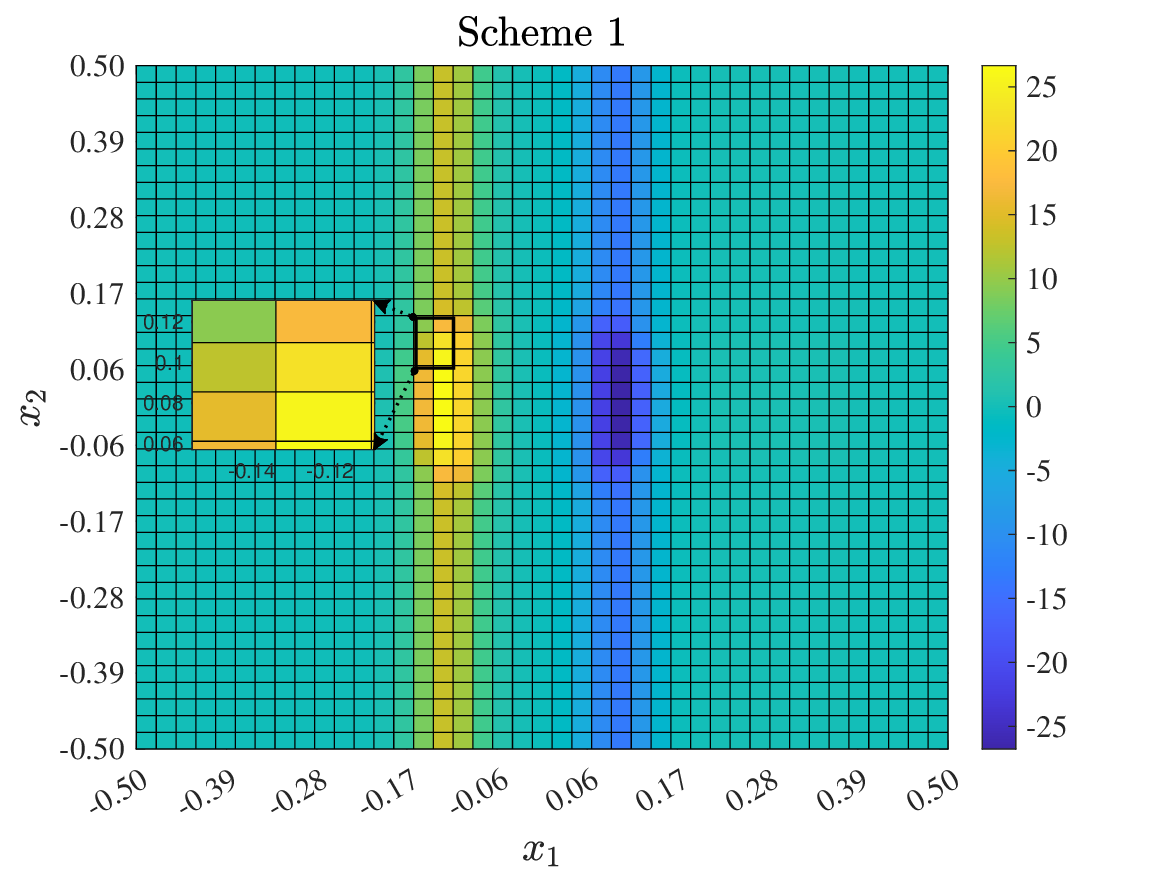} 
    \end{minipage}
    \begin{minipage}[b]{0.49\textwidth}
        \centering
        \includegraphics[width=\textwidth]{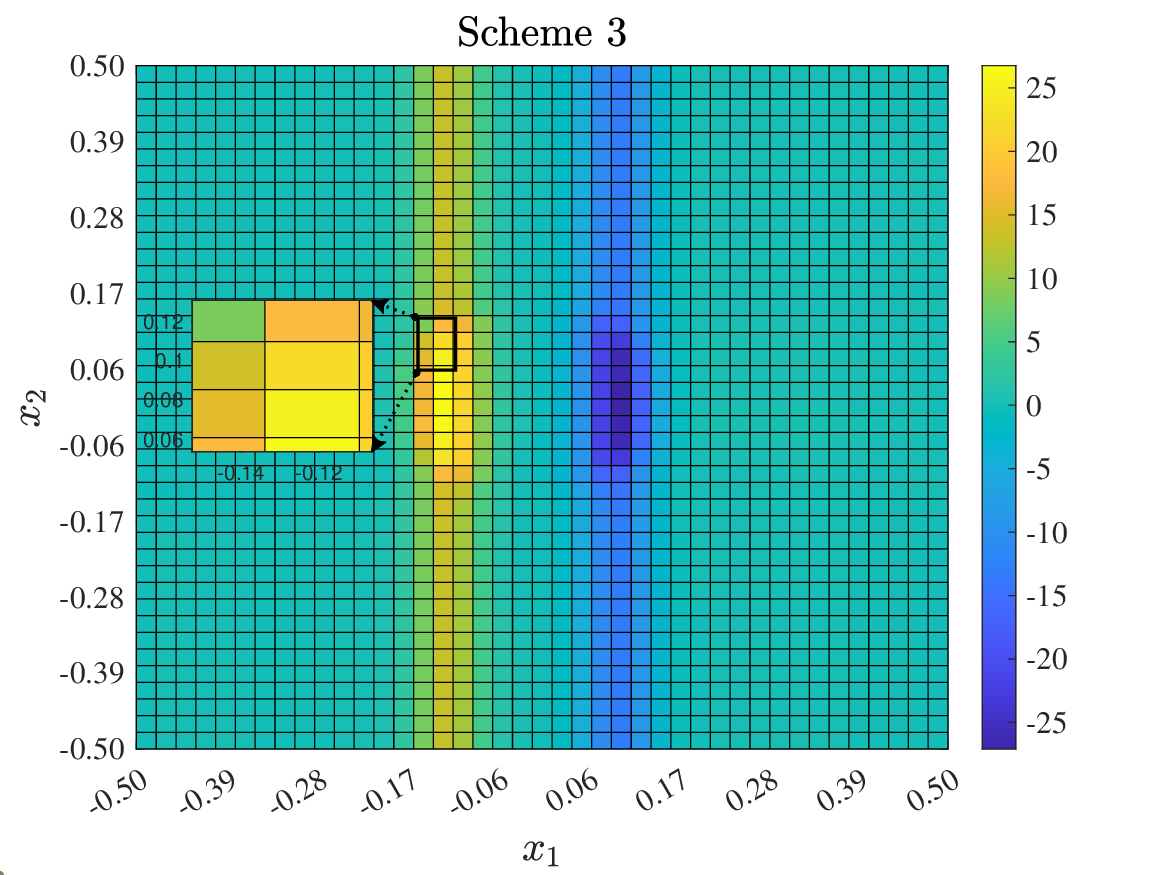}
    \end{minipage}
\caption{Problems III Stage 1. Moment $\partial_{x_2} m_{22}^{\text{NN}}$ with different schemes at final time $T=0.1$.}
    \label{fig:test3_stage1_m21_v3}
\end{figure}

\noindent \textbf{Stage 1. } We first study Stage 1 to compare different schemes with different time steps $T$.
In Figure \ref{fig:test3_stage1_ref_v2}, we plot the reference moments $\partial_{x_1} m_{21}$ and $\partial_{x_2} m_{22}$ solutions at final time $T=0.1$. We plot the approximated solutions $\partial_{x_1} m_{21}^{\text{NN},\textit{stage1}}$  and $\partial_{x_2} m_{22}^{\text{NN},\textit{stage1}}$ generated by different schemes in Figure \ref{fig:test3_stage1_m21_v2} and \ref{fig:test3_stage1_m21_v3}. It can be seen that the preferred scheme 1 can better approximate the reference moments than scheme 3.  To qualitatively demonstrate the performance of scheme 1, we report the relative $\ell^2$ errors generated by different schemes. It can be seen that the relative $\ell^2$ error of scheme 1 is smaller than that of other schemes from  Table \ref{tab:test3_stage1_3}.

\begin{figure}[htbp]
    \centering
    \begin{minipage}[b]{0.49\textwidth}
        \centering
        \includegraphics[width=\textwidth]{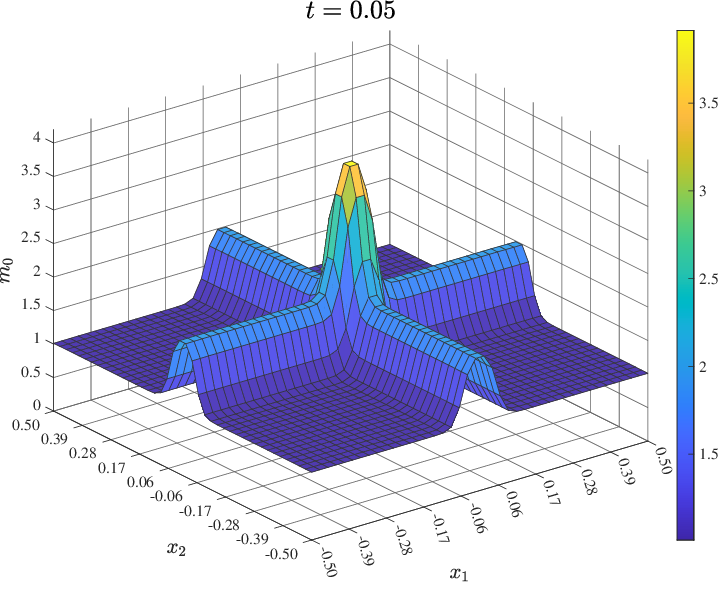} (a)
    \end{minipage}
        \begin{minipage}[b]{0.49\textwidth}
        \centering
        \includegraphics[width=\textwidth]{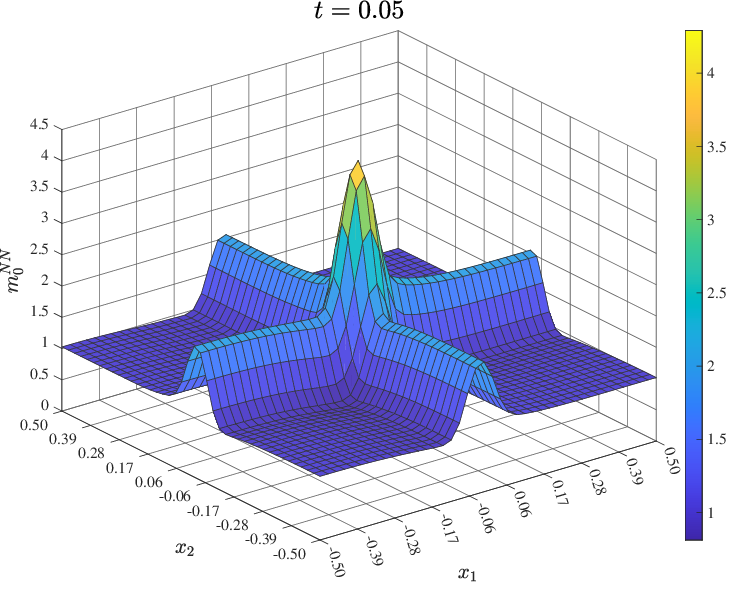} (b)
    \end{minipage}
\caption{Problems III Stage 2.  Moment $m_{0}$ for our proposed method and reference solutions at time step $T=0.05$. }
    \label{fig:test3_stage3_m0}
\end{figure}

\begin{figure}[htbp]
    \centering
    \begin{minipage}[b]{0.49\textwidth}
        \centering
        \includegraphics[width=\textwidth]{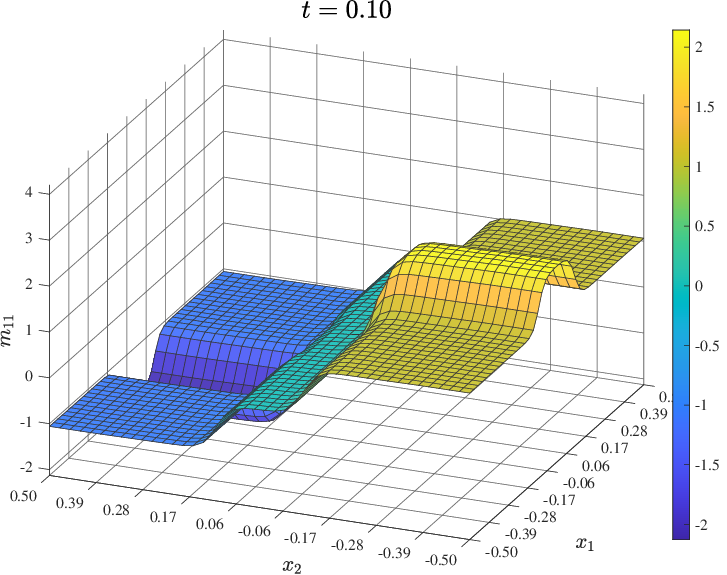} (a)
    \end{minipage}
        \begin{minipage}[b]{0.49\textwidth}
        \centering
        \includegraphics[width=\textwidth]{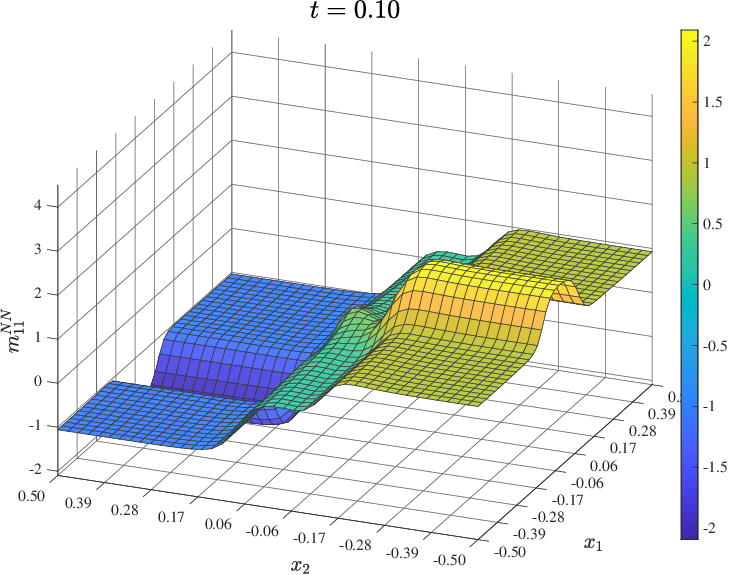} (b)
    \end{minipage}
    \begin{minipage}[b]{0.49\textwidth}
        \centering
        \includegraphics[width=\textwidth]{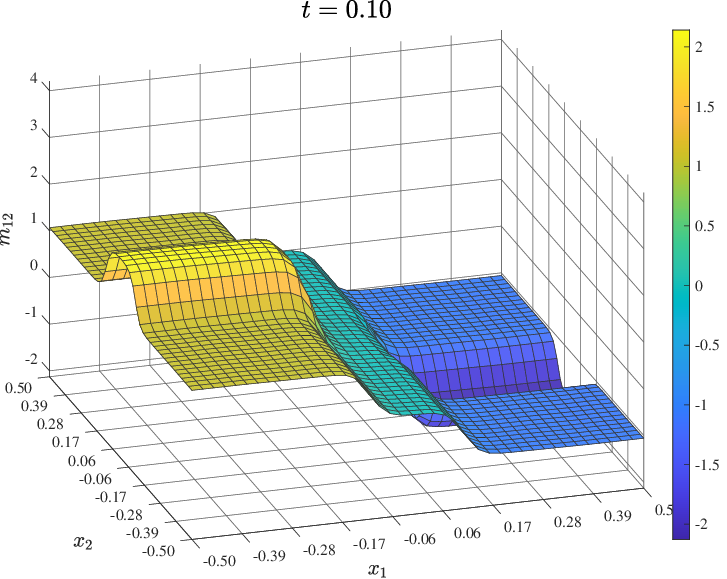} (c)
    \end{minipage}
        \begin{minipage}[b]{0.49\textwidth}
        \centering
        \includegraphics[width=\textwidth]{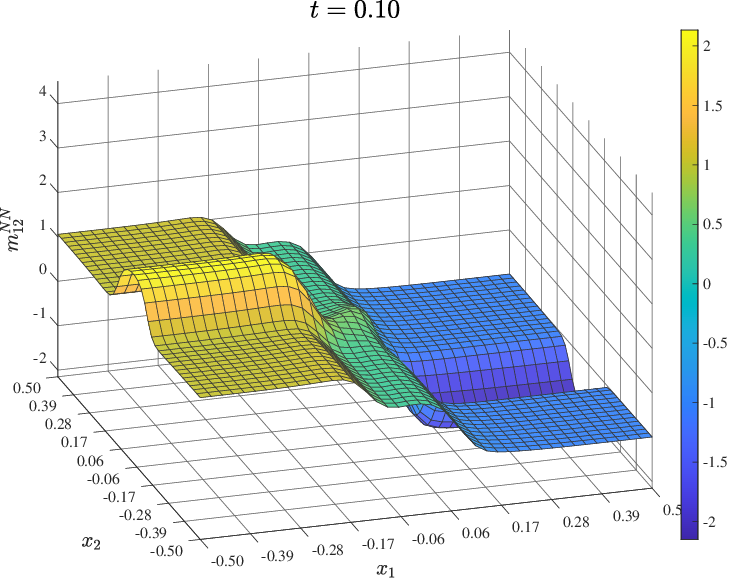} (d)
    \end{minipage}
\caption{Problems III Stage 2.  Moment $m_{11}$ and $m_{12}$ for our proposed method and reference solutions at time step $T=0.10$. 
}
    \label{fig:test3_stage3_m11}
\end{figure}

\begin{table}[htbp]
    \centering
    \begin{small}
            \setlength{\tabcolsep}{8pt} 
            \begin{tabular}{lcccc}
                \toprule
                Error & \textbf{Scheme 1} & Scheme 2 & Scheme 3 & Scheme 4 \\
                \midrule
                $\partial_{x_1}m_{21}$ & $\mathbf{3.120 \times 10^{-2}}$ & $3.175 \times10^{-2}$ & $1.151 \times 10^{-1}$ & $3.200 \times 10^{-2}$ \\
                $\partial_{x_2}m_{22}$ & $\mathbf{3.045 \times 10^{-2}}$ & $3.181 \times 10^{-2}$ & $8.233 \times 10^{-2}$ & $3.632 \times 10^{-2}$ \\
                \bottomrule
            \end{tabular}
            \caption{Problems III.  Relative $\ell^2$ error between the reference moments $\partial_{x_1} m_{21}$ , $\partial_{x_2} m_{22}$ , and the approximated moments $\partial_{x_1} m_{21}^\text{NN}$, $\partial_{x_2} m_{22}^\text{NN}$ at final time $T=0.1$.}
            \label{tab:test3_stage1_3}
    \end{small}
\end{table}
\noindent \textbf{Stage 2.} We further study Stage 2 to compare different moments with different time steps $T$. In Figure \ref{fig:test3_stage3_m0} and \ref{fig:test3_stage3_m11}, we plot reference solutions for $m_{0}$, $m_{11}$, $m_{12}$ and approximated solutions $m_0^{\text{NN},\textit{stage2}}$,\,$m_{11}^{\text{NN},\textit{stage2}}$ and $m_{12}^{\text{NN},\textit{stage2}}$ obtained by by PINNs. We can find from the figures that the moments can match real moments under different time steps well. The detailed results can be found in Table \ref{tab:test3_m0_m1}. Except for the relative $\ell^2$ error, we also study the MSE since moments such as $m_0$ contain large amount points with value equals to zero, which will affect the error results computed by matrices such as relative $\ell^2$ errors.

\begin{table}[htbp]
    \centering
            \setlength{\tabcolsep}{8pt} 
            \begin{tabular}{lcccc}
                \toprule
                Error & $\mathcal{T} = 0.05$, $\ell^2$ &  $\mathcal{T} = 0.10$, $\ell^2$ & $\mathcal{T} = 0.05$, MSE & $\mathcal{T} = 0.10$, MSE \\
                \midrule
                $m_0$ & $1.201 \times 10^{-1}$ & $8.030 \times 10^{-2}$ & $2.406 \times 10^{-2}$ & $1.641 \times 10^{-2}$ \\
                $m_{11}$ & $1.301 \times 10^{-1}$&  $1.510 \times 10^{-1}$&  $2.005 \times 10^{-2}$ & $3.064 \times 10^{-2}$\\
                $m_{12}$ & $1.272 \times 10^{-1}$& $1.461 \times 10^{-1}$ & $1.915 \times 10^{-2}$ & $2.871 \times 10^{-2}$\\
                \bottomrule
            \end{tabular}
            \caption{Problems III Stage 2. Relative $\ell^2$ error between the reference moments $m_0$, $m_{11}$,$m_{12}$ and the approximated moments $m_0^\text{NN}$, $m_{11}^\text{NN}$, $m_{12}^\text{NN}$  at different time steps.}
            \label{tab:test3_m0_m1}
\end{table}

\section{Conclusion and Future Work}\label{sec.4}

We have presented a novel two-stage neural network approach for computing multi-phase solutions to the semiclassical Schrödinger equation. By studying the moment system of the Liouville equation, we developed a method that first learns the relationship between higher- and lower-order moments and their derivatives, then employs physics-informed neural networks (PINNs) to close the system. Theoretical convergence guarantees were established for both the loss function and neural network approximations. Extensive numerical experiments in 1D and 2D demonstrated the method's accuracy across various phase configurations. Future research directions include extending the approach to systems with more phases, higher spatial dimensions, and exploring operator learning techniques for more efficient moment-relation modeling.

\section*{Acknowledgement}
\label{sec:ack}J. W. Jang is supported by the National Research Foundation of Korea (NRF) grants RS-2023-00210484, RS-2023-00219980, and 2021R1A6A1A10042944, and Samsung Electronics Co., Ltd. (IO230412-05903-01).  
J.Y. Lee was supported by Institute for Information \& Communications Technology Planning \& Evaluation (IITP) through the Korea government (MSIT) under Grant No. 2021-0-01341 (Artificial Intelligence Graduate School Program (Chung-Ang University)).
L.~Liu acknowledges the support by National Key R\&D Program of China (2021YFA1001200), Ministry of Science and Technology in China, Early Career Scheme (24301021) and General Research Fund (14303022 \& 14301423) funded by Research Grants Council of Hong Kong from 2021-2023.

\appendix
\section{Details for 2D Problem}
Multiplying \eqref{Vlasov-eqn2D} by $v^l$ ($l=0,1,\cdots, 2N-1$) and integrating over $v$, one can derive the moment equations. Multiplying \eqref{Vlasov-eqn2D} by $v^0$, we get

\begin{equation*}
    \partial_t m_0 + \partial_{x_1} m_{11} + \partial_{x_2} m_{12} = \int_{v_1} \int_{v_2} (x_1 \partial_{v_1} w + x_2 \partial_{v_2} w) d v_2 d v_1.
\end{equation*}
Multiplying \eqref{Vlasov-eqn2D} by $v_1$, we have
\begin{equation*}
    \partial_t m_{11} + \int_{v_1} \int_{v_2} (v_1^2 \partial_{x_1} w  +  v_1 v_2 \partial_{x_2} w) d v_2 d v_1  -  \int_{v_1} \int_{v_2}  v_1 (x_1 \partial_{v_1} w + x_2 \partial_{v_2} w) d v_2 d v_1=0. 
\end{equation*}
Multiplying \eqref{Vlasov-eqn2D} by $v_2$ we have:
\begin{equation*}
    \partial_t m_{12} + \int_{v_1} \int_{v_2} (v_1 v_2 \partial_{x_1} w  +  v_2^2 \partial_{x_2} w) d v_2 d v_1  -  \int_{v_1} \int_{v_2}  v_2 (x_1 \partial_{v_1} w + x_2 \partial_{v_2} w) d v_2 d v_1=0. 
\end{equation*}
\subsection{Schemes for Test III}
\label{scheme_appen}
For scheme 2, we consider
\begin{equation*}
    \begin{cases}
        \partial_{x_1} m_{21}^{\text{NN}, \textit{stage}_1} = \mathcal{NN}_1(\partial_{x_1} m_0, \partial_{x_2} m_0, \partial_{x_1} m_{11}, \partial_{x_2} m_{11},\partial_{x_1} m_{12},\partial_{x_2} m_{12}),\\
        \partial_{x_2} m_{22}^{\text{NN}, \textit{stage}_1} = \mathcal{NN}_2(\partial_{x_1} m_0, \partial_{x_2} m_0, \partial_{x_1} m_{11}, \partial_{x_2} m_{11},\partial_{x_1} m_{12},\partial_{x_2} m_{12}). 
    \end{cases}
\end{equation*}

For scheme 3, we consider
\begin{equation*}
    \begin{cases}
        \partial_{x_1} m_{21}^{\text{NN}, \textit{stage}_1} = {\boldsymbol{c}}^{\mathcal{NN}_1} \cdot (\partial_{x_1} m_0, \partial_{x_2} m_0, \partial_{x_1} m_{11}, \partial_{x_2} m_{11},\partial_{x_1} m_{12},\partial_{x_2} m_{12}),\\
        \partial_{x_2} m_{22}^{\text{NN}, \textit{stage}_1} = {\boldsymbol{c}}^{\mathcal{NN}_2} \cdot (\partial_{x_1} m_0, \partial_{x_2} m_0, \partial_{x_1} m_{11}, \partial_{x_2} m_{11},\partial_{x_1} m_{12},\partial_{x_2} m_{12}), 
    \end{cases}
\end{equation*}
where $\boldsymbol{c}^{\mathcal{NN}_1} = \left( c_i^{\mathcal{NN}_1} = {\mathcal{NN}_1}(m_0,m_{11},m_{12}) \right)_{i=1}^6$ and $\boldsymbol{c}^{\mathcal{NN}_2} = \left( c_i^{\mathcal{NN}_2} \right)_{i=1}^6
$ are coefficient vectors.

For scheme 4, we consider
\begin{equation*}
    \begin{cases}
        \partial_{x_1} m_{21}^{\text{NN}, \textit{stage}_1} = {\boldsymbol{c}}^{\mathcal{NN}_1} \cdot (m_0, \partial_{x_1} m_0, \partial_{x_2} m_0, m_{11}, \partial_{x_1} m_{11}, \partial_{x_2} m_{11}, m_{12}, \partial_{x_1} m_{12}, \partial_{x_2} m_{12})\\
        \partial_{x_2} m_{22}^{\text{NN}, \textit{stage}_1} = {\boldsymbol{c}}^{\mathcal{NN}_2} \cdot (m_0, \partial_{x_1} m_0, \partial_{x_2} m_0, m_{11}, \partial_{x_1} m_{11}, \partial_{x_2} m_{11}, m_{12}, \partial_{x_1} m_{12}, \partial_{x_2} m_{12}), 
    \end{cases}
\end{equation*}
where $\boldsymbol{c}^{\mathcal{NN}_1} = \left( c_i^{\mathcal{NN}_1} \right)_{i=1}^9$ and $\boldsymbol{c}^{\mathcal{NN}_2} = \left( c_i^{\mathcal{NN}_2} \right)_{i=1}^9
$ are coefficient vectors. In Stage 2, we apply PINNs to system \eqref{test3_system} to derive different moments.

\subsection{Empirical Loss for Test III}
Empirical loss for both stage $1$ and stage $2$ of our proposed scheme can be found as follows. For stage $1$, one has
\small{\begin{equation}
\label{R_Stage1}
\begin{aligned}
\mathcal{R}_{\mathrm{Stage 1}} 
&= \sum_{i=1}^{N_1} \left| \partial_{x_1} m_{21}^{\text{NN}}(t_i, x_1^i,x_2^i) - \partial_{x_1} m_{21}^{Data}(t_i, x_1^i,x_2^i) \right|^2 \\
& + \sum_{i=1}^{N_1} \left| \partial_{x_2} m_{22}^{\text{NN}}(t_i, x_1^i,x_2^i) - \partial_{x_2} m_{22}^{Data}(t_i, x_1^i,x_2^i) \right|^2.
\end{aligned}
\end{equation}}
The empirical loss for stage $2$ is as follows:
\small{\begin{equation}
\label{R_Stage2}
\begin{aligned}
&\mathcal{R}_{\mathrm{Stage 2}} \\&= \frac{1}{N_1} \sum_{i=1}^{N_1} \Big| \partial_t m_0^{\text{NN},\textit{stage2}}(t_i, x_1^i, x_2^i) + \partial_{x_1} m_{11}^{\text{NN},\textit{stage2}}(t_i, x_1^i, x_2^i) + \partial_{x_2} m_{12}^{\text{NN},\textit{stage2}}(t_i, x_1^i, x_2^i) \Big|^2 \\
&+ \frac{1}{N_2} \sum_{i=1}^{N_2} \Big| \partial_t m_{11}^{\text{NN},\textit{stage2}}(t_i, x_1^i, x_2^i) + \partial_{x_1} m_{21}^{\text{NN},\textit{stage1}}(t_i, x_1^i, x_2^i) \\&+ \frac{\Phi}{x_1}(t_i, x_1^i, x_2^i) m_{11}^{\text{NN},\textit{stage2}}(t_i, x_1^i, x_2^i) \\
& + \int_{v_1} \int_{v_2} \big( v_1 v_2 \partial_{x_2} w(t_i, x_1^i, x_2^i, v_1, v_2) + v_1 x_2^i \partial_{v_2} w(t_i, x_1^i, x_2^i, v_1, v_2) \big) d v_2 d v_1 \Big|^2 \\
&+ \frac{1}{N_3} \sum_{i=1}^{N_3} \Big| \partial_t m_{12}^{\text{NN},\textit{stage2}}(t_i, x_1^i, x_2^i) + \partial_{x_2} m_{22}^{\text{NN},\textit{stage1}}(t_i, x_1^i, x_2^i) \\&- \frac{\Phi}{x_2}(t_i, x_1^i, x_2^i) m_{12}^{\text{NN},\textit{stage2}}(t_i, x_1^i, x_2^i) \\
& + \int_{v_1} \int_{v_2} \big( v_1 v_2 \partial_{x_1} w(t_i, x_1^i, x_2^i, v_1, v_2) - v_2 x_1^i \partial_{v_1} w(t_i, x_1^i, x_2^i, v_1, v_2) \big) d v_2 d v_1 \Big|^2 \\
&+ \frac{\lambda_1}{N_4} \sum_{i=1}^{N_4} \Big| \mathcal{B} \big( m_0^{\text{NN},\textit{stage2}}(t_i, x_1^i, x_2^i) - m_0(t_i, x_1^i, x_2^i) \big) \Big|^2 \\
&+ \frac{\lambda_2}{N_4} \sum_{i=1}^{N_4} \Big| \mathcal{B} \big( m_{11}^{\text{NN},\textit{stage2}}(t_i, x_1^i, x_2^i) - m_{11}(t_i, x_1^i, x_2^i) \big) \Big|^2 \\
&+ \frac{\lambda_3}{N_4} \sum_{i=1}^{N_4} \Big| \mathcal{B} \big( m_{12}^{\text{NN},\textit{stage2}}(t_i, x_1^i, x_2^i) - m_{12}(t_i, x_1^i, x_2^i) \big) \Big|^2 \\
&+ \frac{\lambda_4}{N_5} \sum_{i=1}^{N_5} \Big| m_0^{\text{NN},\textit{stage2}}(0, x_1^i, x_2^i) - m_0(0, x_1^i, x_2^i) \Big|^2 \\
&+ \frac{\lambda_5}{N_5} \sum_{i=1}^{N_5} \Big| m_{11}^{\text{NN},\textit{stage2}}(0, x_1^i, x_2^i) - m_{11}(0, x_1^i, x_2^i) \Big|^2\\&+ \frac{\lambda_6}{N_5} \sum_{i=1}^{N_5} \Big| m_{12}^{\text{NN},\textit{stage2}}(0, x_1^i, x_2^i) - m_{12}(0, x_1^i, x_2^i) \Big|^2. 
\end{aligned}
\end{equation}}
Here $N_1$, $N_2$, $N_3$, $N_4$, and $N_5$ are the number of sample points of $\mathcal{T} \times \Omega$, $\mathcal{T} \times \Omega$, $\mathcal{T} \times \Omega$, $\mathcal{T} \times \partial \Omega$, and $\Omega$, respectively. For spatial points $(x_1^i, x_2^i)$, we select interior points evenly on $[-0.5, 0.5] \times [-0.5, 0.5]$, and $\mathcal{B}$ stands for the boundary condition. For temporal points $t_i$, interior points are evenly picked in $[0, 0.1]$. The tensor product grid for the collocation points is used, and one sets the penalty parameters in \eqref{R_Stage2} to be $(\lambda_1, \lambda_2, \lambda_3, \lambda_4, \lambda_5, \lambda_6) = (1, 1, 1, 1, 1, 1)$.

\bibliographystyle{siam}
\bibliography{moment.bib}
\end{document}